\documentclass{amsart}
\usepackage{amssymb, amscd}
\numberwithin{equation}{section}

\def\AA{{\mathbb A}}
 
\def\CC{{\mathbb C}}  
\def\DD{{\mathbb D}}

\def\GG{{\mathbb G}}

\def\LL{{\mathbb L}}
\def\NN{{\mathbb N}}
\def\PP{{\mathbb P}}
\def\QQ{{\mathbb Q}} 
 
\def\ZZ{{\mathbb Z}}

\def\g{\gamma}
\def\hs{{\rm h}}
\def\HS{{\rm HS}}
\def\hn{{\rm hn}}
\def\irr{{\rm irr}} 
 
\def\orb{{\rm orb}}
 
\def\reg{{\rm reg}} 
\def\p{{\partial}}

\def\sing{{\rm sing}}
\def\top{{\rm top}}

\def\wt{{\rm wt}}
\def\ra{\rangle}
\def\la{\langle}
\def\bs{\backslash}
\def\add{{\rm add}}

\def\abar{\overline{a}}
\def\bbar{\overline{b}}
\def\lambar{\overline{\lambda}}
\def\lam{\lambda}

\def\Ccal{{\mathcal C}}

\def\Ecal{{\mathcal E}} 

\def\Ical{{\mathcal I}} 
\def\Jcal{{\mathcal J}} 
\def\Khat{\hat{K}}

\def\Lcal{{\mathcal L}}

\def\Mhat{\hat{M}}

\def\mcal{{\mathfrak m}}

\def\Ocal{{\mathcal O}}
\def\Scal{{\mathcal S}}

\def\Vcal{{\mathcal V}}
\def\Xcal{{\mathcal X}}
\def\Ycal{{\mathcal Y}}

\def\mutilde{\tilde{\mu}}
\def\mubold{{\bf\mu}}
\def\muboldhat{{\hat{\bf\mu}}}

\def\la{\langle}
\def\ra{\rangle}

\newcommand\ac{\operatorname{ac}}

\newcommand\aut{\operatorname{Aut}}
\newcommand\conjclass{\operatorname{conj}}
\newcommand\Der{\operatorname{Der}}

\newcommand\Gr{\operatorname{Gr}}
\newcommand\Hom{\operatorname{Hom}}

\newcommand\ord{\operatorname{ord}}

\newcommand\rk{\operatorname{rk}} 

\newcommand\spe{\operatorname{sp}}
\newcommand\spec{\operatorname{Spec}}
\newcommand\supp{\operatorname{supp}}

\newcommand\Tr{\operatorname{Tr}}

\newtheorem{theorem}{Theorem}[section]
\newtheorem{lemma}[theorem]{Lemma}
\newtheorem{keylemma}[theorem]{Key lemma}
\newtheorem{proposition}[theorem]{Proposition}
\newtheorem{corollary}[theorem]{Corollary}

\theoremstyle{definition}
\newtheorem{definition}[theorem]{Definition}

\theoremstyle{remark} 
\newtheorem{remark}[theorem]{Remark}

\newtheorem{question}[theorem]{Question}

\newtheorem{conventions}[theorem]{Conventions}

\title{Motivic measures}
\author{Eduard Looijenga}

\address{Mathematisch Instituut\\
Universiteit Utrecht\\
PO Box 80.010, NL-3508 TA Utrecht\\
Nederland}
\email{looijeng@math.uu.nl}

\begin{document}
\maketitle

\maketitle

\section{Introduction}
An $n$-jet of an arc in an algebraic variety is a one parameter Taylor 
series of length $n$ in that variety. To be precise, if the variety $X$ is
defined over the algebraically closed field $k$, then it is a
$k[[t]]/(t^{n+1})$-valued point of $X$. The set of such
$n$-jets are the closed points of a variety $\Lcal_n(X)$ also defined over
$k$ and the arc space of $X$, $\Lcal (X)$, is the projective limit of
these. Probably Nash \cite{nash} was the first to study arc spaces 
in a systematic fashion (the paper in question was written in 1968).
He concentrated on arcs based at a given point of $X$ and observed
that to each irreducible component of this `provariety' there corresponds in an
injective manner an irreducible component of the preimage of this
point in any resolution of $X$. He asked the (still unanswered) question how to
identify these components on a given resolution. The renewed
interest in arc spaces has a different origin, however. 
Batyrev \cite{11} proved that two connected projective complex 
manifolds with trivial canonical bundle which are birationally equivalent
must have the same Betti numbers. This he showed by first lifting the
data to a situation over a discrete valution ring with finite residue
field and then exploiting a $p$-adic integration technique. (Such a
$p$-adic integration approach to problems in complex algebraic geometry had also
been used by Denef and Loeser \cite{23} in their work on
topological zeta functions attached to singular points of complex
varieties.) When Kontsevich learned of Batyrev's result he saw how this
proof could be made to work in a complex setting using arc spaces. 
The new proof also gave more: equality 
of Hodge numbers, and even an isomorphism of Hodge structures with
rational coefficients. The underlying technique, now going under the name 
of {\it motivic integration}, has led to an avalanche of applications.
These include new (so-called {\it stringy}) invariants of singularities, a
complex analogue of the Igusa zeta function, a motivic version of the
Thom-Sebastiani property and the motivic McKay correspondence. Some of
these were covered in a recent talk by Reid \cite{R} in this seminar. 

The idea is simple if we keep in mind an analogous, more classical
situation. Consider the case of a complete discrete valuation ring $(R,m)$
with finite residue field $F$. There is a Haar measure on the Boolean
algebra consisting of the cosets of powers of $m$ that takes the value $1$
on $R$ (so it is also a probability measure). This induces one on a suitable
Boolean
algebra of subsets of the set of $R$-valued points of any scheme that is
flat of pure dimension and of finite type over $\spec(R)$. Associated to
this measure
is a function that essentially counts the number of `points' in each
reduction modulo $m^k$: the {\it Igusa zeta function}, introduced by Weil,
and intensively studied by Igusa, Denef and Loeser  
(and reported on by Denef in this seminar \cite{22}).
A missing case was that of equal characteristic zero: $\Ocal =k[[t]]$,
$k\supset\QQ$. The proposal of Kontsevich is to give
$\Ocal$ a  measure that takes values in a Grothendieck ring of 
$k$-varieties in which the class of the affine line, $\LL$, is invertible:
the value on the
ideal $(t^n)$ is then simply $\LL^{-n}$ (or $\LL^{1-n}$, which is 
sometimes more convenient). If $\Xcal$ is a suitable $\Ocal$-scheme,
then we obtain a measure on the set of sections as
before, but now with values in this Grothendieck ring. The corresponding
zeta function is a very fine bookkeeping device, for it does its counting
in a ring that is huge. There is no a priori reason to restrict
to the case of equal characteristic, for Kontsevich's idea makes sense for
any complete discrete valuation ring. Indeed, with little extra effort the
material in Sections 2,3 and 9 can be generalized to that context.

This report concerns mainly work of Denef and Loeser. Some 
of their results are presented here somewhat differently, and this is why 
more proofs are provided than one perhaps expects of the write up of a
seminar talk.  References to the sources are in general given
after the section titles, rather than in the statements of theorems.

\smallskip
I thank Jan Denef for inviting me for a short visit to Leuven
to discuss the material exposed here. I am also indebted to Maxim Kontsevich
and especially to Jan Denef for comments on previous versions, from 
which this text has greatly benefitted (though remaining errors are my responsability only). This applies in particular to the 
motivic Thom-Sebastiani theorem
and a word of explanation is in order here. In the original version 
I had introduced (albeit somewhat implicitly) a binary operator on a certain
Grothendieck ring of motives, called here quasi-convolution.  
Quasi-convolution is almost associative,
but not quite, and since I thought this to be a serious defect, I
passed to the universal associative quotient. 
But in a recent overview, Denef and Loeser 
\cite{30} noted that there is no need for this: the 
property one wants (which is another than
associativity) holds already without passing that to quotient. As
this no longer  justifies its introduction, 
I thought it best to take advantage of their
observation and rewrite things accordingly.

\section{The arc space and its measure \cite{25}, \cite{K}}
Throughout the talk we fix a complete discrete valuation ring $\Ocal$ 
whose residue field $k$ is assumed to be algebraically closed and
of characteristic zero. The spectrum of $\Ocal$ is denoted $\DD$ with
generic point $\DD^\times$ and closed point $o$. A uniformizing parameter is
often denoted by $t$ so that $\Ocal =k[[t]]$. The assumption that $k$ be 
algebraically closed is for convenience only: in most situations this
restriction is unnecessary or can be avoided. 

The symbol $\NN$ stands for the set of nonnegative integers.

\medskip
\subsection*{The Grothendieck ring of varieties}
Consider the {\it Grothendieck ring} $K_0(\Vcal_k)$ of reduced $k$-varieties: this
is the abelian group generated by the isomorphism classes of
such varieties, subject to the relations $[X-Y]=[X]-[Y]$, where $Y$ is a
closed in $X$. The product over $k$ turns it into a ring. Note that if we restrict
ourselves to smooth varieties we get the same ring: the reason is that every
$k$-variety $X$ admits a stratification (i.e., a filtration by closed
subschemes $X=X^0\supset X^1\supset\cdots \supset
X^{d+1}=\emptyset$ such that $X^k-X^{k+1}$ is smooth) and that any two
such admit a common refinement. The latter property implies
that $[X]:=\sum_k [X^k-X^{k+1}]$ is unambiguously defined.
In fact, $K_0(\Vcal_k)$ is generated by
the classes of complete nonsingular varieties, for any 
smooth variety $U$ admits a completion $\overline{U}$ by
adding a normal crossing divisor and then $[U]=\sum (-1)^i[\overline{U}^i]$,
where $\overline{U}^i$ stands for the normalization of the codimension $i$
skeleton  of the resulting stratification. W\l odarczyk's weak
factorization theorem (in the form of the main theorem of \cite{wlo}) can be used to show that relations of the following simple type suffice: if $X$ is smooth projective and 
$\tilde X\to X$ is obtained by blowing up a smooth closed subvariety
$Y\subset X$ with exceptional divisor $\tilde Y$, then 
$[\tilde X]-[\tilde Y]=[X]-[Y]$.

We denote the class of the affine line $\AA^1$ by $\LL$ and we write $M_k$
for the  localization $K_0(\Vcal_k)[\LL^{-1}]$. 
Recall that a subset of a variety $X$ is called {\it constructible}
if it is a finite union of (locally closed) subvarieties. Any
constructible subset $C$ of $X$ defines an element $[C]\in M_k$. The
constructible subsets of $X$ form a Boolean algebra and so we
obtain in a tautological manner a  $M_k$-valued measure $\mu_X$
defined on this Boolean algebra. More generally, a morphism $f:Y\to X$
defines on that same algebra an $M_k$-valued measure $f_*\mu_Y$:
assign to a constructible subset of $X$ its preimage in $Y$.

The ring $M_k$ is interesting, big, and hard to grasp. 
Fortunately, there are several characteristics of $M_k$ 
(i.e., ring homomorphisms from $M_k$ to a ring) that are
well understood. We describe some of these in  decreasing order of
complexity under the assumption that $k$ is a subfield of $\CC$. The 
first example is the Grothendieck ring $K_0(\HS)$ of
the category of Hodge structures. A {\it Hodge structure} consists of a 
finite dimensional $\QQ$-vector space $H$, a finite bigrading 
$H\otimes\CC=\oplus_{p,q\in\ZZ}
H^{p,q}$ such that $H^{p,q}$ is the complex conjugate of $H^{q,p}$ and
each {\it weight summand}, $\oplus_{p+q=m}H^{p,q}$, is defined over
$\QQ$. There are evident notions of tensor product and morphism of Hodge
structures so that we get an abelian category $\HS$ with tensor product. 
The Grothendieck construction produces a group $K_0(\HS)$,
elements of which are representable as a formal difference of 
Hodge structures $[H]-[H']$ and $[H]=[H']$ if and only if $H$ and $H'$ are
isomorphic. The tensor product makes it a ring. 

For every complex variety $X$, the cohomology with compact supports, 
$H^r_c(X;\QQ)$, comes with a natural finite increasing filtration
$W_\bullet H^r_c(X;\QQ)$, {\it the weight filtration}, such that the
associated graded $\Gr_\bullet^WH^r_c(X;\QQ)$ 
underlies a Hodge structure having $\Gr_m^WH^r_c(X;\QQ)$ as weight $m$
summand. We assign to $X$ the {\it Hodge characteristic}\footnote{As all 
our characteristics are compactly supported we omit
the otherwise desirable subscript $c$ from the notation.}  
\[
\chi_\hs(X):=\sum_r (-1)^r[H^r_c(X;\QQ)]\in K_0(\HS )
\]
If $Y\subset X$ is closed subvariety, then the exact sequence
\[
\dots\to H^r_c(X-Y)\to H^r_c(X)\to H^r_c(Y)\to H^{r+1}_c(X-Y)\to \dots
\]
is compatible in a strong sense with the Hodge data. This implies
the additivity property
$\chi_\hs(X)=\chi_\hs(X-Y)+\chi_\hs(Y)$. For the affine line
$\AA^1$, $H^r_c(\AA^1;\QQ)$ is nonzero only for $r=2$; the cohomology
group
$H^2_c(\AA^1;\QQ)$ is one-dimensional and of 
type $(1,1)$. So $\chi_\hs(\AA^1)$ (usually denoted as $\QQ(-1)$) is
invertible. It follows that $\chi_\hs$ factorizes over $M_k$.
If we only care for dimensions, then we compose with the ring homomorphism
$K_0(\HS)\to \ZZ[u,u^{-1},v, v^{-1}]$,  $[H]\mapsto \sum_{p,q}
\dim(H^{p,q})u^pv^q$, to get the {\it Hodge number characteristic}
$\chi_\hn :M_k\to \ZZ[u,u^{-1},v, v^{-1}]$. 
It takes $\LL$ to $uv$. The {\it weight characteristic}
$\chi_\wt :M_k\to \ZZ[w,w^{-1}]$ is obtained if we go
further down along the map $\ZZ[u,u^{-1},v, v^{-1}]\to \ZZ[w,w^{-1}]$ 
that sends both $u$ and $v$ to $w$. Evaluating the latter at
$w=1$ gives the ordinary\footnote{A complex algebraic variety can be compactified within its homotopy type by giving it a topological boundary that is stratifyable into strata of odd dimension. This boundary has zero Euler characteristic, hence
the compactly supported Euler characteristic of the variety is its ordinary
Euler characteristic.} Euler characteristic $\chi_\top :M_k\to \ZZ$.

In the spirit of this discussion is the following question raised by  Kapranov \cite{kap}:

\begin{question}
Let $X$ be a  variety over $k$. If $\sigma_n(X)\in M_k$ denotes the class
of its $n$th symmetric power, is then
\[
Z_X(T):=1+\sum_{n=1}^\infty \sigma_n(X)T^n\in M_k[[T]]
\]
a rational function in the sense that it determines an element in a suitable
localization of $M_k[T]$? (Since the logarithmic derivative $Z'/Z$ defines an additive map $M_k\to M_k[[T]]$, we may restrict ourselves 
here to the case of a smooth variety.) Does it satisfy a 
functional equation when $X$ is smooth and complete? Kapranov shows that the answer to both questions is yes in case $\dim (X)\le 1$.
\end{question} 

\subsection*{A measure on the space of sections} Let us call a 
{\it $\DD$-variety} a separated reduced scheme that is flat and
of finite type over $\DD$ and whose closed fiber is reduced. 
Given a $\DD$-variety $\Xcal /\DD$ with closed fiber $X$, then 
the set of its sections up to order $n$, $\Xcal_n$, is the set of
closed points of a $k$-variety (also denoted $\Xcal_n$) 
naturally associated to $\Xcal$. It is obtained from $\Xcal$ modulo
$\mcal^{n+1}$ essentially by Weil restriction of scalars \cite{greenb1}.
So $\Xcal_0=X$. The set $\Xcal_\infty$ of sections of
$\Xcal\to\DD$ is
the projective limit of these and is therefore the set of closed points
of a provariety. 
If $\Xcal/\DD$ is of the form $X\times\DD\to\DD$, with $X$ a $k$-variety,
then we are dealing with the space of {\it $n$-jets} (of curves) on $X$
and the {\it arc space} of $X$, here denoted by $\Lcal_n(X)$ resp.\
$\Lcal (X)$.

For $m\ge n$ we have a forgetful morphism $\pi_n^m: \Xcal_m\to\Xcal_n$.
(When $n=0$, we shall often  write $\pi^m_X$, $\pi_X$ instead of  $\pi_0^m$,
$\pi_0$.) A fiber of $\pi_n^{n+1}$ lies in an affine space over the 
Zariski tangent space of the base point. 
In case $X$ is smooth, it is in fact an affine space 
over the tangent space of the base point: $\pi_n^{n+1}$ has then the
structure of a torsor over the tangent bundle.
A theorem of Greenberg \cite{greenb2} 
asserts that there exists a constant $c$ such that the image of 
$\pi_n$ equals the image of $\pi_n^{cn}$. So $\pi_n(\Xcal_\infty)$ 
is constructible.

The goal is to define a measure on an interesting algebra of
subsets of $\Xcal_\infty$ in such a way that its direct image under
$\pi_X$ is the tautological measure $\mu_X$ when $X$ is smooth.
(This will lead us to deviate from the
definition of Denef-Loeser and Batyrev by a factor $\LL^d$ and to adopt
the one used in \cite{R} instead.) 
For this we assume that $\Xcal$ is of pure relative dimension $d$ and we 
say that a subset $A$ of $\Xcal_\infty$ is {\it stable} if for some
$n\in\NN$ we have 
\begin{enumerate}
\item[---] $\pi_n(A)$ is constructible in $\Xcal_n$ and
$A=\pi_n^{-1}\pi_n(A)$,
\item[---] for all $m\ge n$ the projection $\pi_{m+1}(A)\to
\pi_m(A)$ 
is a piecewise trivial fibration (that is, trivial relative to a 
decomposition into subvarieties) with fiber an affine space of dimension
$d$.
\end{enumerate}
The second condition is of course superfluous in case $\Xcal /\DD$ is
smooth. It is clear that $\dim \pi_m(A)-md$ is independent of
the choice of $m\ge n$; we call this
the {\it (virtual) dimension} $\dim A$ of $A$. The same is true for the class
$[\pi_m(A)]\LL^{-md}\in M_k$; we denote that
class by $\mutilde_{\Xcal}(A)$. The collection of stable
subsets of $\Xcal$ is a Boolean ring (i.e., is
closed under finite union and difference) on which $\mutilde_{\Xcal}$ 
defines a finite additive measure. A theorem of Denef-Loeser 
(see Theorem \ref{noshift}) ensures that there are plenty of stable sets.

In order to extend the measure to a bigger collection of interesting
subsets of $\Xcal_\infty$ we need to complete $M_k$. Given $m\in\ZZ$, let
$F_mM_k$ be the subgroup of $M_k$ spanned by the $[Z]\LL^{-r}$ with
$\dim Z\le m+r$. This is a filtration of $M_k$ as a ring: 
$F_mM_k. F_nM_k\subset F_{m+n}M_k$.
So the separated completion of $M_k$ with respect to this filtration,
\[
\Mhat_k :=\lim_{\leftarrow} M_k/F_mM_k\quad \text{($m\to -\infty$ in
this limit),}
\]
to which we will refer as the {\it dimensional completion}, is also a
ring. The kernel of the natural map $M_k\to \Mhat_k$ is $\cap_m F_mM_k$,
of course. It is not known whether this is zero\footnote{This issue 
is avoided if we work with the adic completion 
$\ZZ((L^{-1}))\otimes_{\ZZ [L]}K_0(\Vcal_k)$ instead, 
but in practice this is too small. Nevertheless, it seems that 
in all applications we are dealing with elements lying in the 
localization $\QQ (L)\otimes_{\ZZ [L]}K_0(\Vcal_k)$.}.
In case $k\subset \CC$, the Hodge characteristic extends to this
completion:  
\[
\chi_\hs :\Mhat_k\to \Khat_0(\HS).
\]
Here $\Khat_0(\HS)$ is defined in a similar way as $\Mhat_k$ with
`dimension' replaced by `weight'. The assertion follows from the fact that
the weights in the compactly supported cohomology of a variety of
dimension $d$ are $\le 2d$. Likewise we can extend the characteristics
counting Hodge numbers or weight numbers (with values Laurent power series
in the reciprocals of their variables). This does not apply to
the Euler characteristic, but in many cases of interest the weight characteristic gives a rational function in $w$ that has no pole at $w=1$. Its value there is then a good substitute.

We will be mostly concerned with the composite of $\mutilde_{\Xcal}$
and  the completion map, for it is this measure that we
shall extend. We call this the {\it motivic measure} on $\Xcal$ and
denote it by $\mu_{\Xcal}$. Let us say that a subset
$A\subset\Xcal_\infty$
is  {\it measurable} if for every (negative) integer $m$
there exist a stable subset $A_m\subset\Xcal_\infty$ and a sequence 
$(C_i\subset\Xcal_\infty)_{i=0}^\infty$ of stable subsets such that the
symmetric difference $A\Delta A_m$ is contained in $\cup_{i\in\NN} C_i$
with $\dim C_i<m$ for all $i$ and $\dim C_i\to -\infty$, for
$i\to\infty$.

\begin{proposition}\label{measurable}
The measurable subsets of $\Xcal_\infty$ make up a Boolean
subring and $\mu_{\Xcal}$ extends as a measure to this ring by
\[
\mu_{\Xcal} (A):= \lim_{m\to -\infty} \mu_{\Xcal}(A_m).
\]
In particular, the above limit exists in $\Mhat_k$
and its value only depends on $A$. 
\end{proposition}

The proof is based on  

\begin{lemma}\label{finitecover}
Let $\Xcal /\DD$ be of pure dimension and $A\subset\Xcal_\infty$ a
stable subset. If $\Ccal =\{C_i\}_{i=1}^\infty$ is a countable covering of
$A$ by
stable subsets with $\dim C_i\to -\infty$ as $i\to\infty$, then $A$ is
covered by a finite subcollection of $\Ccal$.
\end{lemma}
\begin{proof}
Let $n\in\NN$ be such that $A=\pi_n^{-1}\pi_n(A)$. Suppose that $A$ is not
covered by a finite subcollection of $\Ccal$. Choose 
$k\in\NN$ such that $\dim C_i<-(n+2)d$ for $i > k$ and let 
$u_{n+1}\in \pi_{n+1}(A\setminus \cup_{i\le k} C_i)$. We have
$\pi_{n+1}^{-1}u_{n+1}\subset A$. This set is not
covered by a finite subcollection of $\Ccal$, for clearly 
$\pi_{n+1}^{-1}(u_{n+1})$ is not covered by $\{ C_i\}_{i\le k}$
and for $i>k$, $C_i\cap \pi_{n+1}^{-1}(u)$ is of
positive codimension in $\pi_{n+1}^{-1}(u)$. 

With induction we find a sequence
$\{ u_m\in\Lcal_m(X)\}_{m>n}$ so that for all $m>n$ $u_{m+1}$ lies over
$u_m$ and  $\pi^{-1}(u_m)$ is not covered by a finite subcollection of
$\Ccal$. The sequence defines an element $u\in\Xcal$. Since
$\pi_n(u)\in\pi_n
(A)$, we have $u\in A$ and so $u\in C_i$ for some $i$. But if $C_i$ is
stable at level $m>n$, then $\pi_m^{-1}(u_m)\subset C_i$, which
contradicts a defining property of $u_m$.
\end{proof}

For $k=\CC$, the condition $\lim_{i\to\infty}\dim C_i= -\infty$
is unnecessary, for we may then use the Baire property of $\CC$ instead \cite{12}. 

\begin{proof}[Proof of \ref{measurable}]
Suppose we have another solution $A\Delta A'_m \subset
\cup_{i\in\NN} C'_i$ with $A'_m$ and $C'_i$ stable, $\dim(C'_i)<m$ 
for all $i$ and $\dim C'_i\to -\infty$ as $i\to\infty$. 
It is enough to prove that the dimension of the stable set $A_m\Delta
A'_m$ is $<m$. Since $A_m\Delta A'_m\subset \cup_{i\in\NN} (C_i\cup
C'_i)$, Lemma \ref{finitecover} applies and we find that 
$A_m\Delta A'_m\subset \cup_{i\le N} (C_i\cup C'_i)$ for some $N$. 
Since every term has dimension $<m$, this is also true for 
$A_m\Delta A'_m$.
\end{proof}

So a countable union of stable sets $A=\cup_{n\in\NN}
A_n$ with $\lim_{n\to\infty}\dim A_n=-\infty $ is measurable and 
$\mu_{\Xcal} (A)=\lim_{n\to\infty} \mu_{\Xcal} (\cup_{k\le n}
A_k)$.

\begin{remark}
Given a $\DD$-variety $\Xcal$, then for 
any $d\in\NN$ there is a $d$-measure $\mu^d_{\Xcal}$ 
that induces $\mu_{\Ycal}$ on $\Ycal_\infty$ for any 
$\DD$-subvariety $\Ycal$ of pure dimension $d$. We expect this measure to 
extend to a much bigger collection of subsets of $\Xcal$ so that if 
$f: \Xcal\to\Scal$ is a dominant $\DD$-morphism of 
pure relative dimension $d$, then every fiber of 
$f_*:\Xcal_\infty\to\Scal_\infty$ is $\mu^d_{\Xcal}$-measurable.
\end{remark}

Here is a sample of the results of Denef and Loeser on
the rationality of Poincar\'e series \cite{25}.

\begin{theorem}
Let $X$ be a $k$-variety.
Then  $\sum_{n=0}^\infty \mu_{\Xcal} (\pi_n\Lcal(X))T^n\in M_k[[T]]$ 
is a rational expression in $T$ with each factor 
in the denominator of the form $1-\LL^aT^b$ where $a\in\ZZ$ 
and $b$ is a positive integer. 
\end{theorem}

We will not discuss its proof, since this theorem is not used in what 
follows.
Denef and Loeser derive this by means of Kontsevich's transformation rule 
discussed below, which is applied to a suitable projective resolution 
$\Xcal$, and a theorem about semialgebraic sets, due to Pas \cite{pas}. 
It is likely that this theorem still holds for the space
of sections of any $\DD$-variety.

\section{The transformation rule \cite{K}, \cite{25}, \cite{27}} 
We describe two results that are at the basis of the theory. 
The proofs are relegated to Section \ref{proofs}.

\begin{proposition}\label{constructible}
For a $\DD$-variety $\Xcal /\DD$ of pure dimension, the preimage of any 
constructible subset under $\pi_n: \Xcal_\infty\to\Xcal_n$ is measurable. 
In particular, $\Xcal_\infty$ is measurable. If $\Ycal\subset\Xcal$ is 
nowhere dense, then $\Ycal_\infty$ is of measure zero. 
\end{proposition}

For $\Xcal/\DD$ of pure relative dimension we have the notion of an 
integrable function $\Phi: \Xcal_\infty\to \Mhat_k$: this requires the 
fibers of $\Phi$ to be 
measurable and the sum $\sum_a \mu_{\Xcal}(\Phi^{-1}(a))a$ to converge,
i.e., there are at most countably many nonzero terms 
$(\mu_{\Xcal}(\Phi^{-1}(a_i))a_i)_{i\in\NN}$ and we have
$\mu_{\Xcal}(\Phi^{-1}(a_i))a_i\in F_{m_i}\Mhat_k$ with $\lim_{i\to\infty}
m_i=-\infty$. The motivic integral of $\Phi$ is then by definition the 
value of this series:
\[
\int \Phi\, d\mu_{\Xcal} =\sum_i \mu_{\Xcal}(\Phi^{-1}(a_i))a_i .
\]
We have a similar notion for maps with values in topological 
$\Mhat_k$-modules. An important example arises from an ideal 
$\Ical\subset\Ocal_X$: 
such an ideal defines a function $\ord_\Ical:\Xcal_\infty\to \NN\cup\{\infty\}$ by assigning to $\g\in\Xcal_\infty$ the multiplicity of $\g^*\Ical$. The condition $\ord_\Ical\g =n$ only depends on the $n$-jet of $\g$ and 
this defines a constructible subset $C_n\subset\Xcal_n$. Hence the fibers 
of $\ord_\Ical$ are measurable. We shall see that the function 
\[
\LL^{-\ord_\Ical}: \Xcal_\infty\to \Mhat_k
\]
is integrable.  

There is a beautiful transformation rule for motivic integrals under 
modifications. Let $H:\Ycal\to \Xcal$ be a morphism of $\DD$-varieties of 
pure  dimension $d$. We define the {\it Jacobian ideal} 
$\Jcal_H\subset\Ocal_\Ycal$ of 
$H$ as $0$th Fitting ideal of $\Omega_{\Ycal/\Xcal}$. This has the nice 
property that its formation commutes with base change. 
The following theorem generalizes an unpublished theorem of Kontsevich. 

\begin{theorem}\label{transforule} 
Let $H:\Ycal\to\Xcal$ be a $\DD$-morphism of pure 
dimensional $\DD$-varieties with $\Ycal/
\DD$ smooth. If $A$ is a measurable subset of 
$\Ycal_\infty$ with $H\big|_A$ injective, then $HA$ is measurable and 
$\mu_\Xcal(HA)=\int_A \LL^{-\ord_{\Jcal_H}}\, d\mu_\Ycal$. 
\end{theorem} 

\section{The basic formula \cite{25}} 

\subsection*{A relative Grothendieck ring} 
It is convenient to be able to work in a relative setting. 
Given a variety $S$, denote by $K_0(\Vcal_S)$ the 
Grothendieck ring of $S$-varieties and by $M_S$ its localization with 
respect $\LL$. The ring $M_S$ can be dimensionally completed as usual. 
Notice that an element of $M_S$ defines a $M_k$-valued measure on 
on the Boolean algebra of constructible subsets of $S$. Often measures 
are naturally represented this way. For instance, the preceding shows 
that for all $n\in\NN$, the direct image of $\mu_\Xcal$ on
$\Xcal_n$ is given by an element $\mu_{\Xcal ,n}\in\Mhat_{\Xcal_n}$. (Notice that $\mu_{\Xcal ,n}$ is then the direct image of $\mu_{\Xcal ,n+1}$.)

A morphism $f: S'\to S$ induces a
ring homomorphism $f^*: M_S\to M_{S'}$. This makes $M_{S'}$ a $M_S$-module. We also have a direct image $f_*:M_{S'}\to M_S$ that is a homomorphism 
of $M_S$-modules. Notice that $f$ itself defines an element $[f]\in M_S$; 
this is also the image of $1\in M_{S'}$ under $f_*$. 

There are corresponding characteristics. For instance, the ordinary 
Euler characteristic $\chi_\top$ becomes a ring homomorphism 
from $M_S$ to the Grothendieck ring of constructible $\QQ$-vector spaces 
on $S$. This ring is generated by direct images of irreducible local 
systems of $\QQ$-vector spaces over smooth irreducible 
subvarieties $Z$ of $S$. (A
better choice is to take the intersection cohomology sheaf in $S$ of this 
local system along $Z$; this has the advantage that it only depends on the
generic point of $Z$.) 

Similarly, the Hodge characteristic $\chi_\hs$ takes values in a
ring $K_0(\HS_S)$ that is generated by variations of
Hodge structures over a smooth subvariety of $S$. 
The homomorphisms $f^*$ and $f_*$ persist on this level: $f: S'\to S$
induces homomorphisms  $f^* : K_0(\HS_{S})\to K_0(\HS_{S'})$ and
$f_* : K_0(\HS_{S'})\to K_0(\HS_{S})$.

\subsection*{The basic computation} 
A case of interest is when the base variety is $(\NN\times\GG_m)^r$. This 
fails to be finite type, but that is of no consequence and we identify 
$\Mhat_{(\NN\times\GG_m)^r}$ with $\Mhat_{\GG_m^r}[[T_1,\dots ,T_r]]$ in 
the obvious way. 

We use a uniformizing parameter of $\Ocal$ to define 
\[
\ac :\Lcal (\AA^1)-\{ 0\}\to\NN\times\GG_m,
\]
by assigning to $\gamma$ its order $\ord (\g)$ 
resp.\ the first nonzero coefficient of $\g$ ($\ac$ stands 
for {\it angular component}). Integration along $\ac$ sends a
$\Mhat_k$-valued measure on $\Lcal (\AA^1)$ to an element of 
$\Mhat_{\GG_m} [[T]]$. 
The prime example is when this measure is given by a regular 
function $f:\Xcal\to\AA^1$ on a $\DD$-variety $\Xcal$ 
of pure relative dimension: this induces a map 
$f_*:\Xcal_\infty\to\Lcal (\AA^1)$ and we then define 
\begin{equation*}
\begin{CD}
\ac_f: \Xcal_\infty @>{f_*}>>\Lcal(\AA^1) @>{\ac}>> \NN\times\GG_m ,
\end{CD}
\end{equation*} 
so that $[\ac_f] \in \Mhat_{\GG_m}[[T]]$. 
More generally, given a morphism 
$f=(f_1,\dots ,f_r): \Xcal\to\AA^r$, we abbreviate 
\[
\ac_{X,f}:=(\pi_X,\ac_{f_1},\dots \ac_{f_r}): \Xcal_\infty\to X\times
(\NN\times \GG_m)^r.
\]
So $[\ac_{X,f}]\in \Mhat_{X\times\GG_m^r}[[T_1,\dots ,T_r]]$. 

\begin{conventions}\label{resconvent}
If $E$ is a simple normal crossing hypersurface on a smooth 
$k$-variety $Y$, then we adhere to the following 
notation throughout the talk: $(E_i)_{i\in \irr (E)}$ denotes the 
collection of 
irreducible components of $E$ (so these are all smooth by assumption) and 
for any subset $I\subset\irr (E)$, $E^\circ_I$ stands for the locus of 
$p\in\tilde X$ with $p\in E_i$ if and only if  $i\in I$. (With this 
convention, $E^\circ_\emptyset=Y-E$.) We denote the complement of the 
zero section of the normal 
bundle of $E_i$ by $U_{E_i}$ (so this is a $\GG_m$-bundle over $E_i$) 
and $U_I$ designates the fiber product of the bundles $U_{E_i}|
E^\circ_I$, $i\in I$ (a $\GG_m^I$-bundle whose total space has the same 
dimension as $Y$). 

If $\Ecal$ is a simple normal crossing hypersurface on a 
$\DD$-variety $\Ycal/\DD$ with $\Ycal$ smooth, then we shall always 
assume that its union with 
the closed fiber $Y$ has also normal crossings. The notational conventions 
are as above to the extent that restriction or intersection 
with $Y$ is indicated by switching from calligraphic to 
roman font (e.g., $E_i=\Ecal_i\cap Y$). If $Y$ is smooth, then we may 
identify $\irr (E)$ with a subset of $\irr (\Ecal)$. (An equality if 
$\Ecal$ has no component in $Y$.)
\end{conventions}

The following proposition accounts for many of the rationality assertions 
in \cite{25}.

\begin{proposition}\label{basic}
Let $\Xcal/\DD$ be a $\DD$-variety of pure relative dimension and 
$H :\Ycal\to\Xcal$ a resolution of 
singularities. Let $\Ecal$ be a simple normal crossing 
hypersurface on $\Ycal$ that has no irreducible component in the closed 
fiber $Y$. Assume that the Jacobian ideal $\Jcal_H$ of $H$ is 
principal and has divisor $\sum_i (\nu_i-1)\Ecal_i$ (so $\nu_i\ge 1$). 
Let for $\rho =1,\dots ,r$, $f_\rho: \Xcal\to\AA^1$ be a regular function 
such that $f_\rho H$ has zero divisor $\sum_i N_{i,\rho} \Ecal_i$ and 
put $N_i:=(N_{i,1},\cdots N_{i,r})\in\NN^r$, $i\in \irr(E)$. Then  
\[
[\ac_{X,f}] =\sum_{I\subset \irr (E)} [U_I/ X\times \GG_m^r] 
\prod_{i\in I} (\LL^{\nu_i}T^{-N_i}-1)^{-1} \text{ in } 
\Mhat_{X\times\GG_m^r}[[T_1,\dots ,T_r]],
\]
where $U_I\to X\times \GG_m^r$ has first component projection onto 
$E^\circ_I\subset X$ followed by the restriction of $H$ and second 
component induced by $fH$. 
\end{proposition}
\begin{proof}
Given $m\in \NN^{\irr (E)}$, consider the set $\Ycal(m)$ of 
$\g\in\Ycal_\infty$ with order $m_i$ along $\Ecal_i$. So for 
$\g\in \Ycal (m)$ we have $\ord_{\Jcal_H}(\g)=\sum_i
m_i(\nu_i-1)$ and $\ord_{f_\rho H} (\g)=\sum_i m_iN_{\rho ,i}$. 
If $\supp (m)\subset \irr (E)$ is the support of $m$, then we  have a 
natural projection $e_m: \Ycal(m)\to U_{\supp (m)}$. 
Its composite with the morphism $U_{\supp (m)}\to X\times\GG_m^r$
is a restriction of $\ac_{X,fH}:=(\pi_X H,\ac_{f_1H},\dots
,\ac_{f_rH}):\Ycal_\infty
\to X\times (\NN\times \GG_m)^r$ with $\NN^r$-component 
$\sum_i m_iN_i$. In other words, 
\[
[\ac_{X,fH }\big|_{\Ycal(m)}]=
[U_{\supp (m)}/ X\times\GG_m^r] \LL^{-\sum_i m_i}T^{\sum_i m_iN_i}. 
\]
So the  transformation formula \ref{transforule} yields 
\begin{align*}
[\ac_{X,f}] &=\sum_{m\in \NN^{\irr (E)}} 
[U_{\supp (m)}/ X\times\GG_m^r] \prod_{i\in\supp (m)}\Big(
\LL^{-m_i-m_i(\nu_i-1)}T^{m_iN_i}\Big)\\
&=\sum_{I\subset \irr (E)} [U_I/ X\times \GG_m^r]
\prod_{i\in I} \frac{\LL^{-\nu_i}T^{N_i}}{1-\LL^{-\nu_i}T^{N_i}}. 
\end{align*}
\end{proof}

If we drop the assumption that $\Ecal$ has no irreducible component in $Y$,
then the above formula must be somewhat modified: now each irreducible 
component of $Y$ contributes with an expression of the above form 
times a monomial in $\LL^{-1}$ and $T_1,\dots ,T_r$.

\begin{corollary}\label{basiccor}
In the situation of \ref{basic}, the class of $(\pi_X,\ord_f):\Xcal_\infty
\to X\times\NN^r$ in $\Mhat_X[[T_1,\dots ,T_r]]$ equals 
\[
\sum_{I\subset \irr (E)} [E^\circ_I/ X]
\prod_{i\in I} \frac{\LL-1}{\LL^{\nu_i}T^{-N_i}-1}.
\]
In particular, the direct image of $\mu_{\Xcal}$ on $X$ is 
represented by
\[
\sum_{I\subset \irr (E)} [E_I^\circ/X]
\prod_{i\in I}[\PP^{\nu_i-1}]^{-1}.
\] 
\end{corollary}
\begin{proof} Since $U_I$ is a $\GG_m^I$-bundle over $E^\circ_I$, the 
class of the projection $U_I\to X$ is $(\LL-1)^{|I|}$ times the class of
$E^\circ_I\to X$. 
\end{proof}

This corollary shows that $\Xcal_\infty$ is measurable so that the
measurable subsets of $\Xcal_\infty$ form in fact a Boolean algebra. 
It also implies that the Hodge number characteristic of $\Xcal_\infty$ is an element of 
$\QQ [u,v][(uv)^N-1)^{-1}\, |\,N=1,2,\dots ]$ on which the Euler characteristic  takes the value
$\sum_{I\subset \irr (E)} \chi_\top( E_I^\circ )\prod_{i\in I}\nu_i^{-1}$.

\begin{remark} We can also express the direct image of $\mu_{\Xcal}$
on $X$ in terms of the closed subvarieties $E_I$: if 
$\irr'(E)$ denotes the set of $i\in\irr (E)$ with $\nu_i\ge 2$, then
\[
\sum_{I\subset\irr'(E)}(-\LL)^{|I|}[E_I/X]\prod_{i\in I}
\frac{[\PP^{\nu_i-2}]}{[\PP^{\nu_i-1}]}.
\]
All varieties appearing in this expression are  proper over $X$ and 
nonsingular. So it gives rise to an element of a complex cobordism ring of $X$ 
localized away from the classes of the complex projective varieties. 
This class, and the values that various genera take on it, might
deserve closer study. 
\end{remark}

\section{The motivic nearby fiber \cite{24}, \cite{29}}

\subsection*{An equivariant Grothendieck ring} 
Let $G$ be an affine algebraic group. We consider varieties $X$ with 
{\it good} $G$-action, where `good' means that 
every orbit is contained in an affine open subset. 
For instance, a representation of $G$ on a $k$-vector space $V$  is good. 
For a fixed variety $S$ with  
$G$-action, we define the  Grothendieck group $K_0^G(\Vcal_S)$ as
generated by isomorphism types of $S$-varieties with good $G$-action
modulo the usual equivalence relation (defined by pairs) and the
relation that declares that every finite dimensional representation
$\rho$ of $G$ has the same class as the trivial representation of the 
same degree (i.e., $\LL^{\deg (\rho )}$). 

In case the action on $S$ is trivial, the 
product makes $K_0^G(\Vcal_S)$ a $K_0(\Vcal_S)$-algebra. 
If moreover $G$ is finite abelian, then assigning to a variety $X$
with good $G$-action its  $G$-orbit space $\overline{X}:=G\bs X$ augments this as a
$K_0(\Vcal_S)$-module: 
\[
K_0^G(\Vcal_S)\to K_0(\Vcal_S),\quad  a\mapsto \bar a .
\]
(Not as an algebra, for the orbit space of a product is in general not the 
product of orbit spaces.) That this is well-defined follows from the
lemma below. (We do not know whether this holds for arbitrary finite $G$.)

\begin{lemma}\label{trivial} 
Let be given a representation of a finite abelian group $G$ on a 
$k$-vector space $V$ of finite dimension $n$. Then the class of 
$\overline{V}$ in $K_0(\Vcal_k)$ is $\LL^n$. 
\end{lemma}
\begin{proof}
Let $V=\oplus_{\chi\in\hat G} V_\chi$ be the eigenspace 
decomposition of the $G$-action. Given a subset $I\subset \hat G$, denote 
by $V_I$ the set of vectors in $V$ whose $V_\chi$-component is 
nonzero if and only if $\chi\in I$. We have a natural projection 
$\overline{V_I}\to\prod_{\chi\in I}\PP (V_\chi)$. This has the structure 
of a torus bundle, the torus in question being a quotient of $\GG_m^I$ 
by a finite subgroup. So the class of $\overline{V_I}$ in $M_k$ is 
$(\LL-1)^{|I|}$ times the class of $\prod_{\chi\in I}\PP (V_\chi)$. 
Since $V_I$ has also that structure, the classes of $\overline{V_I}$ and 
$V_I$ in $M_k$ coincide. Hence the same is true for $\overline{V}$ and $V$.
\end{proof}

Similarly we can form $M_S^G:=K_0^G(\Vcal_S)[\LL^{-1}]$ 
and its dimensional completion. The class of
an $S$-variety $Z/S$ with $G$-action in 
$M_S^G$ or $\Mhat_S^G$ is denoted by $[Z/S;G]$.   
If $G$ is abelian and acts trivially on $S$, then we have 
corresponding augmentations taking values in $M_S$ and its completion.

There are corresponding characteristics in case $k\subset\CC$. For 
instance, the ordinary Euler characteristic defines a ring homomorphism 
from $M_k^G$ to the Grothendieck ring $K^G_0(\QQ )$ of finite 
dimensional representations of $G$ over $\QQ$ and more generally, we have 
a ring homomorphism $\chi^G_\top$ from $M_S^G$ to the Grothendieck ring of 
constructible sheaves with $G$-action on $S$, 
$K^G_0(\QQ_S)$. Similarly, there is a Hodge character 
$\chi^G_\hs: M_S^G\to  K_0^G(\HS_S)$. 

\subsection*{The case $G=\muboldhat$} 
We will mostly (but not exclusively) be concerned with the case 
when $G$ is a group of roots of unity. We have 
the Grothendieck ring $M_S^{\muboldhat}$ 
of varieties with a topological action of the 
procyclic group $\muboldhat=\lim_{\leftarrow}\mubold_n$ (such an action
factorizes through a finite quotient $\mubold_n$). 
The inverse automorphism of $\muboldhat$, $\zeta\mapsto\zeta^{-1}$, 
defines an involution $^*$ in $M_S^{\muboldhat}$. 

The group of continuous characters of $\muboldhat$ is naturally isomorphic 
with $\QQ/\ZZ$, with the involution $^*$ acting as multiplication by $-1$; 
the projection $\muboldhat\to\mubold_n$ followed by the 
inclusion $\mubold_n\subset\GG_m$ corresponds to $\frac{1}{n}\pmod{\ZZ}$. 
In other words, $K^{\muboldhat}_0(\CC)\cong \ZZ[e^\alpha \, |\, \alpha\in
\QQ/\ZZ]$. For every positive integer $n$ there is a rational irreducible 
representation 
$\chi_n$ of $\mubold_n$, namely the field $\QQ (\mubold_n)$, regarded as 
$\QQ$-vector space. These make up an additive basis of 
$K^{\muboldhat}_0(\QQ)$. The image of $\chi_n$ in $\ZZ[\QQ/\ZZ]$ is 
$\sum_{(k,n)=1} e^{k/n}$, which allows us to regard 
$K^{\muboldhat}_0(\QQ)$ 
as a subring of $\ZZ[e^\alpha \, |\, \alpha\in\QQ/\ZZ]$. 

The so-called mapping torus construction gives rise to an
$M_k$-linear map  
\[
M_k^{\muboldhat}\to M_{\GG_m} 
\]
with the property that composition with the direct image homomorphism 
$M_{\GG_m}\to M_k$ is $(\LL-1)$ times the augmentation $M_k^{\muboldhat}\to M_k$. It is defined as follows. If $X$ is a variety with good
$\mu_n$, then its {\it mapping torus} is  the \'etale locally
trivial fibration  $\GG_m\times^{\mu_n}X\to\GG_m$ whose total space is the orbit
space of the $\mu_n$-action on $\GG_m\times X$ defined by 
$\zeta (\lambda ,x)=(\lambda\zeta^{-1},\zeta x)$ and for which the
projection is induced by $(\lambda ,x)\mapsto\lambda^n$. Notice that the
fiber over $1\in\GG_m$ can be identified with $X$ and that the monodromy is
given by the  action of $\mubold_n$ on $X$. The projection on the
second factor induces a morphism $\GG_m\times^{\mu_n}X\to\GG_m\to
\overline{X}$ that has the structure of a piecewise $\GG_m$-bundle.
So the image of $\GG_m\times^{\mu_n}X$ in $M_k$ is
$(\LL-1)[\overline{X}]$. If $m=kn$ is a positive multiple
of $n$ and we let $\mu_m$ act on $X$ via $\mu_m\to\mu_n$, then 
$(\lambda ,x)\mapsto (\lambda^k,x)$ identifies the two fibrations
over $\GG_m$ and so we have a map as asserted. This generalizes 
at once to the case
where we have a base variety with trivial $\muboldhat$-action.

\subsection*{$\aut (\DD)$-equivariance}
The automorphism group $\aut (\DD)$ can be identified with the
group of formal power series $k[[t]]$ with nonzero constant term where the 
group law is given by substitution. It acts on the arc space of any 
$k$-variety by composition: $h(\g):=\g h^{-1}$. 
If the variety  is of pure dimension, then this action is free 
outside negligible subset. Clearly, a morphism of $k$-varieties
induces an $\aut (\DD)$-equivariant map 
between their arc spaces. Since we end up with more than just a $\aut (\DD)$-invariant measure
on an arc space, it is worthwhile to explicate this structure by
means of a  definition. If $\DD_n$ denotes the subscheme of $\DD$
defined by the  ideal $(t^{n+1})$, then $\aut (\DD_n)$ (which has the same underlying variety as the group of units of $k[[t]]/(t^{n+1})$) acts
naturally on $\Lcal_n(X)$. For $n\ge 1$, the kernel of 
$\aut (\DD_{n+1})\to \aut (\DD_n)$ can be identified with $\GG_a$.
Its action is trivial on $(\pi^{n+1}_1)^{-1}(0)$ and free on the complement $(\pi^{n+1}_1)^{-1}(TX-\{ 0\})$. By choosing a constructible section of the latter we lift the direct image homomorphism $(\pi^{n+1}_n)_*$ to a map  
\[
M_{\Lcal_{n+1}(X)}^{\aut (\DD_{n+1})}\to
M_{\Lcal_{n}(X)}^{\aut (\DD_{n})}.
\]
The result is easily seen to be independent of this choice.

\begin{definition}
An {\it equivariant motivic measure} on
$\Lcal(X)$ is a collection $\lam =(\Lambda_n\in
\Mhat_{\Lcal_n(X)}^{\aut (\DD_n)})_{n=1}^\infty$, so that
$\Lambda_n$ is the direct image of $\Lambda_{n+1}$ for all $n$. 
\end{definition}

It is clear that such a collection determines an $\Mhat_k$-valued measure on
the stable subsets. The definition is so devised that the measure $\mu_{\Lcal (X)}$ constructed earlier comes from an equivariant motivic measure.

This notion is of particular interest when the variety in question
is a smooth curve $C$ and we are given a closed point $o\in C$.
An $\aut (\DD)$-orbit in $\Lcal (C,o)$ is given
by a positive integer $n$ that may also take the value $\infty$. 
If $n$ is finite, then this orbit projects onto the set of nonzero 
elements of $(\pi^n_{n-1})^{-1}(0)\cong T_{C,o}^{\otimes n}$. 
The group $\aut (\DD_n)$ acts on the latter orbit through $\aut (\DD_1)\cong \GG_m$ with $\mu_n\subset\GG_m$ as isotropy group. 
So the value of $\lam$ on a fiber over 
$T_{C,o}^{\otimes n}-\{ 0\}$ is naturally an element $\lam_n$ of
$M_k^{\mu_n}$. We call the generating series 
\[
\lam (T):=\sum_{n=1}^\infty \lambda_nT^n
\]
the {\it zeta function} of $\lambda$. 
It is not hard to verify that this series determines $\lam$ completely.
This is particularly so if we view $\lam$ as a $\Mhat_k$-valued measure
on $\Lcal (C,o)$. For instance, its value
on the preimage in $\Lcal (C,o)$ of a constructible subset $A$ of 
$\Lcal_m(C,o)$ consisting of order $n$-arcs (with $n\le m$) is 
$\LL^{n-m}[A]\lambar_n$. 
Notice that the series $\sum_{n=1}^\infty (\LL-1)\lambar_n$
converges to the full integral of $\lambda$.

\subsection*{A motivic zeta function}
Given  a pure dimensional variety $X$ and a flat morphism $X\to \AA^1$, let
$X_0:=f^{-1}(0)$ and denote by $f$ the restriction $(X,X_0)\to (\AA^1,0)$.
Then the direct image of $\mu_{\Lcal (X,X_0)}$ (regarded as an equivariant measure) on $X_0\times\Lcal (\AA^1,0)$ is then also 
equivariant. We will (perhaps somewhat ambiguously) 
refer to this measure as the direct image of the 
motivic measure of $\Lcal (X,X_0)$ on $X_0\times\Lcal(\AA^1,0)$. 
Its zeta function is denoted by 
\[
S(f)=\sum_{n=1}^\infty S_n(f)T^n \in \Mhat_{X_0}^\muboldhat[[T]].
\]
We now assume that $X$ is smooth and connected. The smoothness of $X$ ensures that the preimage in $\Lcal (X,X_0)$ of a stable subset of 
$\Lcal(\AA^1,0)$ of level $n$ is stable of level $n$, so that $S_n(f)$ already is defined as an element of $M_{X_0}$ (but we shall not bring out the
distinction in our notation). The series $S(f)$ can be computed from
an embedded resolution of the zero  set of $f$, $H:Y\to X$ of $X$, 
as in \ref{resconvent}. We assume here that 
the preimage $E$ of $X_0$ is a simple normal crossing hypersurface 
that contains the exceptional set.
Let $m$ be a positive integer that is divided by all the  coefficients $N_i$
of the divisor $(f)$ on the irreducible components of $E$.
If we make a base change of $\tilde f:=fH$ over 
the $m$th power map $\AA^1\to\AA^1$ and normalize, then we get a 
$\mubold_m$-covering $\tilde Y\to Y$. Let $\tilde E_I^\circ$ be a connected
component of the preimage of $E_I^\circ$ in $\tilde Y$. The 
restriction $\tilde E_I^\circ\to E_I^\circ$ is unramified,
and has $\mubold_m$-stabilizer of $\tilde E_I^\circ$ as its Galois group. The 
latter is easily seen to be the subgroup $\mubold_{N_I}$, 
where $N(I):=\gcd\{N_i\,|\, i\in I\}$. This defines 
\[
[\tilde E_I^\circ/ Y;\mu_{N(I)}]\in M_Y^{\mubold_{N_I}}.
\] 
This element lies over $X_0$ if $I$ is nonempty, an assumption we make
from now on. We wish to compare it 
with $U_I(1)\subset U_I$, the fiber over $1$ of the projection 
$U_I\to\GG_m$ induced by $\tilde f$. This projection has weights 
$(N_i)_{i\in I}$ relative to  the $\GG_m^I$-action and so 
$\prod_{i\in I}\mubold_{N_i}\subset \GG_m^I$ preserves 
$U_I(1)$. This finite group contains a monodromy action by 
$\mu_{N(I)}$: write $N(I)=\sum_{i\in I} \alpha_iN_i$ and embed 
$\GG_m$ in $\GG_m^I$ by $t\mapsto (t^{\alpha_i})_{i\in I}$ 
(since the $(\alpha_i)_{i\in I}$ are relatively prime, this is an 
embedding indeed). Notice that the projection $U_I\to\GG_m$ is homogeneous of 
degree $N_I$ relative to the action of this one parameter subgroup. This 
implies that $\mu_{N(I)}\subset\GG_m$ may serve as monodromy group. 
(There are a priori several choices for this action, but they are all 
$E_I^\circ$-isomorphic.) 

\begin{lemma}\label{normalclass} 
In $M_{Y_0}^{\muboldhat}$ we have $[U_I(1)/ Y_0;\mubold_{N(I)}] 
=(\LL-1)^{|I|-1} [\tilde E_I^\circ /Y_0;\mu_{N(I)}]$. 
\end{lemma} 
\begin{proof} 
One verifies that the Stein factorization of the projection $U_I(1)\to
E_I^\circ$ has $\tilde E_I^\circ \to E_I^\circ$ as finite factor with 
$U_I(1)\to \tilde E_I^\circ$ being an algebraic torus bundle
of rank  $|I|-1$. In view of Lemma \ref{trivial} the equivariant class of
the latter is $(\LL-1)^{|I|-1}$ times the equivariant class of the base.
The lemma follows. 
\end{proof}

Much of the work of Denef-Loeser on motivic integration centers around the
following 

\begin{theorem}\label{mainformula}
The following identity holds in $M_{X_0}^{\muboldhat}[[T]]$: 
\[ 
S(f)=\sum_{\emptyset\not= I\subset \irr(E)} (\LL-1)^{|I|-1}
[\tilde E_I^\circ /X_0;\mu_{N(I)}] \prod_{i\in I} 
(\LL^{\nu_i}T^{-N_i}-1)^{-1}.
\]
\end{theorem}
\begin{proof} 
Start with the identity of Proposition \ref{basic} (with $r=1$). 
Omit at both sides the constant terms (on the right this amounts to
summing over nonempty $I$ only), and restrict the resulting identity to
the fiber over $1\in\GG_m$. If we take into account the monodromies and
use Lemma \ref{normalclass}, we get the asserted identity, at least if we
take our coefficients in $\Mhat_{X_0}^\muboldhat$. 
Inspection of the proof shows that this actually holds in 
$M_{X_0}^\muboldhat[[T]]$. 
\end{proof}

So the expression at the righthand side is independent of the resolution,
something that is not at all evident a priori. 
Since it lies in the $M_{X_0}^\muboldhat$-subalgebra of
$M_{X_0}^{\muboldhat}[[T]]$ generated by the fractions 
$(L^\nu T^{-N}-1)^{-1}$ with $\nu, N>0$,  
$S(f)$ has a value at $T=\infty$: 
\[
S(f)\big|_{T=\infty}=-\sum_{\emptyset\not= I\subset \irr(E)}
(1-\LL)^{|I|-1}
[\tilde E_I^\circ/X_0;\mu_{N(I)}] 
\]

\subsection*{Comparison with ordinary monodromy}

The element $-S(f)\big|_{T=\infty}$ has an interpretation in terms of
the nearby cycle sheaf of $f$ as we shall now explain. 

Suppose first that $k=\CC$. Let $\widetilde{X-X_0}\to X-X_0\subset X$ be
the pull-back along $f$ of the universal covering 
$\exp :\CC \to \CC^\times\subset\CC$. Take the full
direct image of the constant sheaf $\QQ_{\widetilde{X-X_0}}$ on $X$ and
restrict to $X_0$: this defines $\psi_f$ as an element of the derived
category of constructible sheaves on $X_0$. Let 
$\sigma :\widetilde{X-X_0}\to \widetilde{X-X_0}$ be a generator of the
covering transformation that induces in $\CC$ translation over $-2\pi\sqrt{-1}$.
This generator has the property that its action in $\psi_f$ is the
monodromy.

Let $H: Y\to X$ be a resolution as in \ref{resconvent}.
In the same way, $\psi_{\tilde f}$ is defined as an element of the
derived category of constructible sheaves on the zero set $Y_0$ of $\tilde
f$. The full direct image of $\psi_{\tilde f}$ on $X_0$ is equal to $\psi_f$.

An elementary calculation shows that the stalk of $\psi_{\tilde f}$ at a
point of
$E_I^\circ$ is the cohomology of $N_I$ copies of a real torus of dimension
$N_I-1$. More precisely, the restriction of $\psi_{\tilde f}$ to
$E_I^\circ$
is naturally representable as the full direct image of the constant sheaf
on $U_I(1)$ (an algebraic torus bundle of dimension $N_I-1$ over 
$\tilde E_I^\circ$) under the projection $U_I(1)\to E_I^\circ$.
We have a canonical isomorphism $H^k(\GG_m^r;\QQ)\cong
H_c^{k+r}(\GG_m^r;\QQ)$ and hence the Euler
characteristic $\sum_k (-1)^k[H^k(\GG_m^r;\QQ)]$ in $K_0(HS)$ is $(-1)^r$
times the Euler characteristic $\sum_k (-1)^k[H^k_c(\GG_m^r;\QQ)]$. 
In other words, it is the value of $\chi_h$ on $(1-\LL)^r$.
Hence, if $Z$ is a subvariety of $E_I^\circ$ with preimage $\tilde Z$ in
$\tilde E_I^\circ$, then $\sum_k (-1)^k[H_c^k(Z;\psi_{\tilde f})]$ is the
value of $\chi_h$
on $(1-\LL)^{|I|-1} [\tilde Z;\mu_{N(I)}]$. This shows that $\psi_{f}$ and
$-S(f)\big|_{T=\infty}$ have the same Hodge characteristic. We therefore put
\[
[\psi_f]:=-S(f)\big|_{T=\infty}=\sum_{\emptyset\not= I\subset
\irr(E)} (1-\LL)^{|I|-1} [\tilde E_I^\circ /X_0;\mu_{N(I)}]. 
\]
We refer to $[\psi_f]$ as the {\it nearby cycle class} of $f$ along $X_0$.
Its component in the augmentation submodule, 
\[
[\phi_f]:= [\psi_f]-\overline{[\psi_f]}\in M^\muboldhat_{X_0},
\]
is by definition the {\it vanishing cycle class} of $f$.

\smallskip
Let $S$ be a variety with trivial $\muboldhat$-action. 
Given a $S$-variety $Z$ with a good topological 
$\muboldhat$-action, then for any 
positive integer $n$ the fixed point locus of $\ker
(\muboldhat\to\mubold_n)$ in $Z$ is a $S$-variety which inherits a good
$\mubold_n$-action. This defines a homomorphism of $M_S$-algebras
\[
\Tr_n : M_S^\muboldhat\to M_S^{\mubold_n}.
\]
If $\sigma\in\muboldhat$ generates a dense subgroup of 
$\muboldhat$, then the fixed point locus of
$\ker (\muboldhat\to\mubold_n)$ is also the fixed point locus of $\sigma^n$. In case $k\subset \CC$, a Lefschetz fixed point formula (applied to a partition of $Z$ by orbit type) implies that $\chi_h[Z^{\sigma^n}]$ equals the trace of $\sigma^n$ in $\chi_h[Z]$.
So we may then think of $Tr_n[Z]$ as the motivic trace of $\sigma^n$. 
This is why the following proposition is a motivic version of a result of
A'Campo \cite{acampo}.

\begin{proposition}[see Denef-Loeser \cite{29}]
The series $S(f)$ and $\sum_{n=1}^\infty \Tr_n[\psi_f]T^n$ in
$M_{X_0}^{\mubold_n}[[T]]$ are congruent modulo $\LL-1$. 
\end{proposition} 
\begin{proof}
The monodromy $\sigma$ acts on $\tilde E^\circ_I$ as a covering
transformation of order $N_I$. So $\sigma^n$ has no fixed point if
$N_I$ does not divide $n$ and is equal to all of $\tilde E^\circ_I$
otherwise. It follows from formula for the nearby cycle class that 
\[
\Tr_n[\psi_f]=\sum_{I\subset\irr
(E),N_I|n }(1-\LL)^{|I|-1} [\tilde E^\circ_i /X_0]. 
\]
So 
\begin{multline*}
\sum_{n=1}^\infty \Tr_n[\psi_f]T^n=
\sum_{n=1}^\infty \sum_{I\subset \irr (E),N_I|n}
(1-\LL)^{|I|-1} [\tilde E^\circ_I /X_0] T^n\\
=\sum_{\emptyset\not= I\subset \irr (E)}
(1-\LL)^{|I|-1} [\tilde E^\circ_I /X_0]\sum_{k\ge 1}T^{kN_I}\\
=\sum_{\emptyset\not= I\subset \irr (E)}
(1-\LL)^{|I|-1} [\tilde E^\circ_I /X_0]\frac{T^{N_I}}{1-T^{N_I}}.
\end{multline*}
If we reduce modulo $(\LL-1)$ only the terms with $I$ a singleton remain.
Theorem \ref{mainformula} shows that this has the same reduction modulo
$(\LL-1)$ as $S(f)$.
\end{proof}

\section{The motivic zeta function of Denef-Loeser \cite{24}}

This function is a motivic analoge of Igusa's local zeta function.
It captures slightly less than the function $S(f)$, but has the virtue
that it is defined in greater generality. 
First we introduce two homomorphisms of Grothendieck rings.

An arrow $M_S^{\mubold_{rn}}\to M_S^{\mubold_n}$ is defined by 
assigning to a variety with good $\mubold_{rn}$-action its orbit space
with respect to the subgroup $\mubold_r\subset\mubold_{rn}$ (with
a residual action of $\mubold_n$). The totality of these arrows forms a
projective system whose limit we denote by $M_S(\muboldhat)$. 
This is not the same as $M_S^{\muboldhat}$, but
there is certainly a natural ring homomorphism
\[
\rho :M_S^{\muboldhat}\to M_S(\muboldhat).
\] 
It is given by assigning to a variety $X$ with good
$\muboldhat$-action, the system
$(X_n)_n$, where $X_n$ is the orbit space of $X$ by the kernel of
$\muboldhat\to \mubold_n$ endowed with the residual action of
$\mubold_n$. 

\smallskip
We next define the {\it Kummer map}
\[
M_{S\times\GG_m}\to M_S(\muboldhat ), \quad [f]\mapsto [f]^{1/\infty}.
\]
Given a $S$-variety $Y$ and a morphism $f:Y\to \GG_m$, then for
every positive integer $n$, let $f^{1/n}: Y(f^{1/n})\to \GG_m$ be the
pull-back of $f$ over the $n$th power map $[n]:\GG_m\to\GG_m$. 
So $Y(f^{1/n})$ is the hypersurface in $\GG_m\times Y$
defined by $f(z)=u^n$. The projection of $Y(f^{1/n})\to Y$ is a
$\mubold_n$-covering and thus defines an element $[f]^{1/n}$ of
$M_S^{\mubold_n}$.  Notice that $Y(f^{1/n})$ is the
orbit space of $Y(f^{1/nr})$ relative to the subgroup
$\mubold_r\subset\mubold_{rn}$. Hence the $[f]^{1/n}$'s define an element 
$[f]^{1/\infty}\in M_S(\muboldhat)$.

The following lemma is a straightforward
exercise.

\begin{lemma}\label{kummerprocess}
The composition of the mapping torus
construction and the Kummer map is equal to $(\LL-1)\rho$.
\end{lemma}

For $\Xcal$ a smooth $\DD$-variety of pure relative dimension $d$, define
the
{\it  Denef-Loeser zeta function} by
\[
I(f):=\LL^{-d}\sum_{n=0}^\infty [\ac_{f,n}]^{1/\infty}\,\LL^{-sn}\in
M_{X}(\muboldhat )[[\LL^{-s}]],
\]
where $\LL^{-s}$ is just a variable with a suggestive notation. Then
Corollary \ref{basiccor} and Lemma \ref{kummerprocess} yield

\begin{theorem}\label{igusathm} 
The following identity holds in $M_X(\muboldhat)[[\LL^{-s}]]$: 
\[ 
I(f)=\LL^{-d}\sum_{I\subset \irr(E)} \rho [\tilde E_I^\circ /
X;\mu_{N(I)}]
\prod_{i\in I} \frac{\LL-1}{\LL^{\nu_i+sN_i}-1}.
\]
\end{theorem}

\subsection*{Putting $\LL=1$} 
Consider the $\ZZ [L,L^{-1}]$-subalgebra
$S$ of $\QQ (L,L^{-s})$ generated by the rational functions
$(L-1)(L^{n+sN}-1)^{-1}$, $n,N\ge 1$.
The spectrum of $S$ contains the  generic point of the exceptional
divisor of the blow up of $(1,1)$ in $\GG_m\times\AA^1$. The
corresponding specialization is the evalation homomorphism 
$S\to \QQ (s)$ which sends $(L-1)(L^{n+sN}-1)^{-1}$ to $(n+sN)^{-1}$. 
According to Theorem \ref{igusathm}, $I(f)$ lies in 
$S\otimes _{\ZZ [L,L^{-1}]} M_X(\muboldhat)$. Evaluation at
$\LL=1$ yields
\[
I(f)\big|_{\LL=1}=\sum_{I\subset \irr(D)} \rho [\tilde
E_I^\circ /X;\mubold_{N(I)}] \prod_{i\in I} \frac{1}{\nu_i+sN_i}\in 
M_X(\muboldhat)/(\LL -1)\otimes_\ZZ\QQ(s). 
\]
This is the motivic incarnation of the topological
zeta function considered earlier by Denef and Loeser in \cite{23}. At
the time the resolution independence of this function was established  using
Theorem \ref{igusa}  below.

\subsection*{Comparison with Igusa's $p$-adic zeta function}
Suppose we are given a complete discrete valuation ring $(R,m)$ of
characteristic zero whose residue field $F=R/m$ has finite cardinality
$q$. Then $R$ contains all the $(q-1)$st roots of unity $\mu_{q-1}$ and 
this group projects isomorphically onto $F^\times$. Let $K$ be the quotient field of $R$. If we choose a uniformizing parameter $\pi\in m-m^2$, then 
then the collection $(\zeta\pi^k)_{\zeta\in\mu_{q-1},k\in\ZZ}$ is a system
of representatives of $K^\times/(1+m)$. Define 
\[
\ac^s :K \to \ZZ[\mu_{q-1}][q^{-s}]
\]
by assigning to $u\in\zeta\pi^k +m^{k+1}$ the value $\zeta q^{-ks}$ and
$0$ to $0$. (Here $q^{-s}$ is just the name of a variable; the righthand 
side can be more canonically understood as the group algebra of $K^\times/(1+m)$.)
There is a natural (additive) Haar measure $\mu$ on the the Boolean ring
of subsets of $K$ generated by the cosets of powers of $m$ that takes 
the value $1$ on $R$. It takes values in $\ZZ [q^{-1}]$.
Given an $f\in R[x_1,\dots ,x_d]$ whose reduction mod $m$ is nonzero, then
its {\it Igusa local zeta function} is defined by
\[
Z(f):=\int_{R^m} \ac^s f (x) d\mu (x),
\]
where $R^m$ is endowed with the product measure. We regard this as an
element of $\QQ [\mubold_{q-1}][[q^{-s}]]$: the coefficient of 
$\zeta q^{-ns}$ is the volume of $f^{-1}(\zeta\pi^n +m^{n+1})$.  
(It is customary to let $s$ be a complex number---the series then converges 
in a right half plane---and to compose with a
complex character $\mubold_{q-1}\to\CC^\times$.) 

Let us write $\Xcal$ for $\spec(R[x_1,\dots ,x_d])$ and regard $f$ as
a morphism  $\Xcal\to\AA^1_R$ over $\spec(R)$.
Suppose we have an embedded resolution $H:\Ycal\to\Xcal$ of the zero locus 
of $f$ over $\spec(R)$ with a simple normal crossing hypersurface $\Ecal$
relative to $\spec(R)$ (so no irreducible component in the closed fiber).
Then we get an embedded resolution of the closed fiber $Y\to X$ with simple
normal crossing divisor $E$.  
Make a base change of $fH:\Ycal\to\AA_R^1$
over the $(q-1)$st power map $[q-1]:\AA^1_R\to\AA^1_R$ and normalize; this
gives a $\mubold_{q-1}$-covering $\tilde\Ycal\to\Ycal$. We now get a covering $\hat E ^\circ_I\to E^\circ_I$ defined over $F$ with 
Galois group $\mu_{N_q(I)}$, where $N_q(I):=\gcd (q-1, (N_i)_{i\in I})$  
over $F$ in much the same way as before. The $\mu_{N_q(I)}$-set 
$\hat E^\circ_I(F)$ determines an element 
\[
\#[\hat E^\circ_I;\mu_{q-1}]\in \QQ [\mu_{N_q(I)}]\subset\QQ [\mu_{q-1}],
\] 
where the last inclusion is defined by the surjection 
$\mu_{q-1}\to\mu_{N_q(I)}$.
Denef proved earlier \cite{21} the following analogue
of \ref{igusathm}:

\begin{theorem}[Denef]\label{igusa}
In this situation we have
\[
Z(f)=q^{-d}\sum_{I\subset \irr(E)}
\# [\hat E^\circ_I;\mu_{q-1}]\prod_{i\in I}\frac{q-1}{q^{\nu_i+sN_i}-1},
\]
where $\nu_i$ and $N_i$ have the usual meaning.
\end{theorem}

As appears from \ref{igusathm}, $Z(f)$
is what we get from the value of $I(f_{\bar K})$ on $X_0(\bar K)$ 
(with $\bar K$ an algebraic closure of $K$) if we 
replace classes in $M_{\bar K}$ by the number of $F$-rational points in
their $F$-counterparts (so that we substitute $q$ for $\LL$) and 
pass from $\muboldhat$ to  $\mubold_{q-1}$. This should be understood on a 
more conceptual level that involves a Grothendieck ring 
$M^{\mubold_{q-1}}_{\spec (R)}$ which specializes to both
$M^{\mubold_{q-1}}_{\spec (\bar K)}$ and 
$\QQ [\mubold_{q-1}][[q^{-s}]]$, and avoids resolution.

\section{Motivic convolution \cite{26}}

\subsection*{Join and quasi-convolution}
Consider the Fermat curve  $J_n$ in $\GG_m^2$ defined
by $u^n+v^n=1$.  Notice that it is invariant under the subgroup
$\mubold_n^2\subset\GG_m^2$. If $d$ is a positive divisor of $n$,
then the $\mubold_d^2$-orbit space
of $J_n$ is $J_{n/d}$. In particular, the $\mubold_n^2$-orbit space
of $J_n$ is $J_1$, an affine line less two points. Given varieties
$X$ and $Y$ with good $\mu_n$-action, then we have the variety with
$\mubold_n\times\mubold_n$-action
\[
J_n(X,Y):= J_n\times^{(\mubold_n\times\mubold_n)} (X\times Y).
\]
(If a group $G$ acts well on varieties $A$ and $B$, then
$A\times ^GB$ stands for quotient of $A\times B$ by the
equivalence relation $(ga,b)\sim (a,gb)$ with $G$ acting well on it
by $g[a,b]:= [ga,b]=[a,gb]$.)
Let $\mu_n$ act on $J_n(X,Y)$ diagonally: $\zeta [(u,v),(x,y)]:= 
[(\zeta u,\zeta v), (x,y)]$. The natural map $J_n(X,Y)\to J_1$ is
\'etale locally trivial. If $Y$ has trivial $\mubold_n$-action, then
$J_n(X,Y)=J_n(X,pt)\times Y$ and the variety $J_n(X,pt)$ can be
identified with $(\GG_m-\{\mubold_n\})\times^{\mubold_n} X$. The
latter has the structure of a piecewise $\GG_m$-bundle over $\overline{X}$
from which a copy of $X$ has been removed. 
Similarly, the natural projection of $\overline{J_n(X,Y)}\to
\overline{X}\times \overline{Y}$ is a piecewise
$\GG_m$-bundle from which a copy of $\overline{X\times Y}$ has been
removed.

The construction is perhaps better understood in terms of the 
fibrations over $\GG_m$ defined by the mapping torus construction. 
Recall that  for a variety $X$ with $\mubold_n$-action, its mapping
torus $\GG_m\times^{\mubold_n}X$ fibers over $\GG_m$ by
$[\lam ,x]\mapsto \lam^n$ with $\{ 1\}\times X$ mapping to the
fiber over $1$. The monodromy is the given $\mubold_n$-action on $X$.
If $Y$ is another variety with $\mubold_n$-action, then the composite
\begin{equation*}
\begin{CD}
(\GG_m\times^{\mubold_n}X)\times (\GG_m\times^{\mubold_n}Y)
@>>> \GG_m\times\GG_m \subset \GG_a\times\GG_a @>{+}>> \GG_a 
\end{CD}
\end{equation*}
is a fibration over $\GG_m$. The fiber over $1\in\GG_a$ is identified
as $J_n(X,Y)$ and the monodromy is the given $\mubold_n$-action on 
$J_n(X,Y)$ defined above. 

Clearly, $J_n(X,Y)\cong J_n(Y,X)$.
If $m$ is a divisor of $n$ and the action of 
$\mubold_n$ on $X$ and $Y$ is through $\mubold_m$, then
$J_m(X,Y)=J_n(X,Y)$. So this induces a binary operation, the \textit{join} 
\[
J: M_k^{\muboldhat}\times M_k^{\muboldhat}\to M_k^{\muboldhat}.
\]
The preceding discussion shows that the join is
commutative and bilinear over $M_k$ and that
(i) $J(a,1)=(\LL-1)\abar -a$ and 
(ii) $\overline{J(a,b)}= (\LL-1)\abar\bbar-\overline{ab}$,
where we recall that $a\in M_k^\muboldhat\mapsto \abar\in M_k$ is the
augmentation defined by `passing to the orbit space'. This suggests
to define another binary operation $*$,  
the \textit{quasi-convolution}, on $M_k^{\muboldhat}$ by:
\[
a*B:= -J(a,b)+(\LL -1)\overline{ab}.
\]
The quasi-convolution is commutative and bilinear over $M_k$, whereas the
properties (i) and (ii)  come down to
\begin{enumerate}
\item[(i)] $1$ is a unit for $*$: $a*1=a$ (and hence
$a*\overline{b}=a\overline{b}$) and 
\item[(ii)] $\overline{a*b}=\overline{ab}$.
\end{enumerate}
Neither the join nor the quasi-convolution is associative, but we do have:
\begin{enumerate}
\item[(iii)] $a*(b*c)-(\LL-1)\overline{a(b*c)}+(\LL-1)^2\overline{a}
\overline{bc}$ is symmetric in $a$, $b$ and $c$,
\end{enumerate}
which shows that the quasi-convolution is associative modulo elements of $M_k$.
This property is seen as follows.
Let $J_n^2$ denotes  the Fermat surface in $\GG_m^3$ defined 
by $u^n+ v^n +w^n=1$ and consider the morphism
\[
J_n\times J_n\to J_n^{(2)},\quad 
((u_1,v_1),(u_2,v_2))\mapsto (u,v,w)=(u_1,v_1u_1,v_1u_2).
\]
This morphism is equivariant with respect to the
action of $\mubold_n$ on $J_n\times J_n$ that is diagonal on the
first factor and trivial on the second and the diagonal action
$\mubold_n$ on $J_n^{(2)}$. It also factorizes over  
the orbit space of $J_n\times J_n$ with respect to the $\mubold_n$
action defined by $\zeta ((u_1,v_1), (u_2,v_2))= ((u_1,\zeta
^{-1}v_1), (\zeta u_2,\zeta  v_2))$. One easily verifies that this
identifies the orbit space for this action with  
in $J_n^{(2)}-K_n$, where $K_n\subset J_n^{(2)}$ is defined by $u^n=1$. 
A choice of an $n$th root $\alpha$ of $-1$, 
identifies $K_n$ with $\mubold_n\times\GG_m\times\mubold_n$ via
$(u,v,w)\mapsto (u,v, \alpha w/v)$.
The $\mubold_n^3$-action on $K_n$ carries in an
obvious manner to $\mubold_n\times\GG_m\times\mubold_n$. 

It follows from these observations that if $X,Y,Z$ are varieties with good 
$\mubold_n$-action, then $J_n^{(2)}\times^{\mubold_n^3}
X\times Y\times Z$ decomposes as a $\mubold_n$-variety into two
pieces that can be identified with $J_n(X,J_n(Y,Z))$ and  
$X\times(\GG_m\times^{\mubold_n}(Y\times Z))$ respectively. 
The factor $\GG_m\times^{\mubold_n}(Y\times Z)$ has the structure of
a $\GG_m$-bundle over $\overline{Y\times Z}$. Passing now to
$M_k^\muboldhat$ we find that
\[
J(a,J(b,c))+ (\LL-1)a\overline{bc}
\]
is symmetric in $a,b,c$ and this is equivalent to property (iii)
above.

Join and quasi-convolution extend to $\Mhat_k^\muboldhat$ and 
admit relative variants.

\subsection*{Formation of the spectrum}
Join and quasi-convolution also descend to the Grothendieck ring  
$K_0^\muboldhat (HS)$ of Hodge structures with $\muboldhat$-action. 
We need:

\begin{lemma}[Shioda-Katsura, \cite{sk}] 
Given $(\alpha ,\beta)\in (\QQ/\ZZ)^2$, then for every common 
denominator $n$ of $\alpha$ and $\beta$, the Hodge type of the 
eigenspace $I_{\alpha ,\beta}$ of $\mu_n\times\mu_n$ in $H^1_c(J_n)$ with
character $(\alpha ,\beta)\in (n^{-1}\ZZ/\ZZ)^2$ is independent of $n$ and
we have $\dim I_{\alpha ,\beta}=1$ for $(\alpha ,\beta)\not=(0,0)$ and  
$\dim I_{0,0}=2$. If $\alpha\in \QQ/\ZZ \mapsto \tilde\alpha\in
[0,1[$  is the obvious section, then  
\[
I_{\alpha,\beta} \text{ is of Hodge type }
\begin{cases}
(0,1) \text{ if } \alpha\not= 0\not=\beta \text{ and }
0<\tilde\alpha +\tilde\beta <1,\\
(1,0)\text{ if } 1<\tilde\alpha +\tilde\beta<2 \text{ and}\\
(0,0) \text{otherwise: }\alpha=0 \text{ or } \beta=0 \text{ or }
\alpha +\beta=0.  
\end{cases}
\]
The only other nonzero group is  $H^2_c(J_n)$, which is isomorphic to 
$\QQ(-1)$ and has trivial 
character $(0,0)$. 
\end{lemma}

\begin{corollary}\label{weightshift}
If $H, H'\in K_0^\muboldhat (HS)$, then
\[
H*H'=H_0\otimes H'_0 + \sum_{\alpha\not= 0} H_\alpha\otimes H'_{-\alpha} (-1)+
\sum_{\alpha +\beta\not= 0} H_\alpha\otimes H'_\beta\otimes I_{\alpha,\beta}.
\]
\end{corollary}

Anderson \cite{anderson} investigated Hodge
structures with $\muboldhat$-action using a notion of a {\it
fractional Hodge structure}. For us such a structure will consist of a
complex vector space $V$ defined over $\QQ$ with a complex 
decomposition  $V=\oplus_{p,q\in\QQ; p+q\in\ZZ} V^{p,q}$ such that $V^{q,p}$ is
the complex conjugate of $V^{p,q}$ and $\oplus_{p+q=n}V^{p,q}$ is defined over $\QQ$ for every $n\in\ZZ$. They form an abelian category $HS(\QQ)$ with
tensor product.  Anderson associates to a Hodge structure $H$ with
$\muboldhat$-action a fractional Hodge structure $\sigma (H)$ whose 
underlying vector space is $H$, leaves the bidegrees on $H_0$
unaltered and increases the bidegrees of $H_\alpha$ by 
$(\tilde\alpha,1-\tilde\alpha)$ if $\alpha\not=0$. We shall refer to
this operation as the \textit{formation of the spectrum}. It
defines an additive functor and hence a homomorphism of groups 
$\spe :K_0^\muboldhat (HS)\to K_0(HS(\QQ))$. This is not a ring homomorphism, 
but Corollary \ref{weightshift} shows that $\spe$ takes 
quasi-convolution to the tensor product:
\[
\spe (H*H')=\spe (H)\otimes\spe (H').
\]

\subsection*{Convolution} 
In what follows we need the (additive) group structure on the 
affine line, so we write $\GG_a$  instead of $\AA^1$. We have a 
bijection $\Lcal(\GG_a,0)\cong\mcal$, defined by assigning to
$\g\in\Lcal(\GG_a,0)$ the pull-back of the standard coordinate on $\GG_a$.

Let $\lam =(\Lambda_n)_n$ and $\lam'= (\Lambda'_n)_n$ be  
equivariant measures on $\Lcal(\GG_a,0)$. Then
$\lam \times \lam':=(\Lambda_n\times\Lambda'_n)_{n=1}^\infty$ defines
a measure on the algebra of stable subsets of $\Lcal(\GG_a,0)^2$ 
(that is, preimages of constructible subset of some truncation
$\Lcal_n(\GG_a,0)^2$). For instance, if $C\subset\Lcal_n(\GG_a,0)^2$  is constructible
and consists of pairs of truncated arcs of fixed order $(k,l)$ (with
$k,l\le n$),  then the value of $\lam\times\lam'$ on the preimage of $C$ in 
$\Lcal(\GG_a,0)^2)$ is $\lam_k\lam'_l[C]\LL^{-2n}$.

The direct image of $\lam \times \lam'$ under the addition morphism
$\add : \GG_a\times\GG_a \to\GG_a$, 
$\lam *\lam':=(\Lcal_n(\add )(\Lambda_n\times\Lambda'_n))_{n=1}^\infty$, 
is an equivariant measure on $\Lcal(\GG_a,0)$, called the 
{\it convolution} of $\lam$ and $\lam'$.

\begin{lemma}\label{convolution}
The zeta function of $\lam *\lam'$ is determined by those of 
$\lam$ and $\lam'$:
\[
(\lam *\lam')_n=-(\lam_n * \lam'_n)+(\LL-1)\sum_{i\le n}
\LL^{i-n}\overline{\lam_i\lam '_i)} 
+ (\LL-1)\sum_{i> n} (\lam_n\overline{\lam'_i}
+\overline{\lam_i}\lam '_n). 
\] 
\end{lemma}
\begin{proof} 
The preimage of $t^n+\mcal^{n+1}$ in
$\mcal\times\mcal$ under $\Lcal (\add)$ decomposes into the following pieces:
$(t^n+\mcal^{n+1})\times \mcal^{n+1}$, $\mcal^{n+1}\times
(t^n+\mcal^{n+1})$ and for $i=1,\dots ,n$ the preimage $\tilde C_{n,i}$ of
the subset $C_{n,i}\subset
((\mcal^i-\mcal^{i+1})/\mcal^{n+1})^2$ of pairs 
$(\alpha_i t^i+\cdots +\alpha_nt^n,\beta_i t^i+\cdots +\beta_nt^n)$
with $\alpha_k+\beta_k=0$ for $k=i,\dots ,n-1$ and $\alpha_n+\beta_n=1$. 
We must evaluate
$\lam\times\lam'$ on each of these (relative to the
diagonal $\mu_n$-action). The first piece gives
$\lam_n\sum_{i>n}(\LL-1)\overline{\lam '_i}$ and the second 
the same expression with $\lam$ and $\lam'$
interchanged. Since $[C_{n,i}]=[(
\mcal^i-\mcal^{i+1})/\mcal^{n+1}]=(\LL-1)\LL^{n-i}$, 
we find that for $i<n$, the value of 
$\lam\times\lam'$ on $\tilde C_{n,i}$ equals 
$(\LL-1)\LL^{i-n}\overline{\lam_i\lam '_i}$
(the action of $\mu_n$ is trivial here).
Notice that $C_{n,n}$ is embedded in 
$(\mcal^n-\mcal^{n+1}/\mcal^{n+1})^2$ as 
$J_1$ in $\GG_m^2$. From the above
discussion one sees that $\lam\times\lam'$ takes on this
set the value $J(\lam_n,\lam'_n)$. If we substute
the defining equation for $*$, the Lemma follows.
\end{proof}

This lemma suggests a notion of a convolution operator for 
series \[
\lam (T)= \sum_{n=1}^\infty \lam_nT^n\in\Mhat_k^\muboldhat [[T]]
\]
with the property that the \textit{mass} $(\LL-1)\sum_{n=1}^\infty \lambar_n$ converges.  

For a $\ZZ [L,L^{-1}]$-module $M$ we set 
\[
M\la T\ra :=M[T][\frac{1}{T^N-L^\nu}\, |\, \nu\in\ZZ, N=1,2,3,\dots ].
\]
Expanding the denominators $(1-T^NL^{-\nu})^{-1}$ in $T$ embeds 
$M\la T\ra$ in $M[[T]]$ and expanding 
$(1-T^{-N}L^\nu)^{-1}$ in $T^{-1}$ embeds $M\la T\ra$ in $M[[T^{-1}]][T]$.

According to Theorem \ref{mainformula}, $S(f)\in\Mhat^\muboldhat \la T\ra $. 

\begin{theorem}[Abstract Thom-Sebastiani property]\label{abstractts}
Let $\lam$ and $\lam$ be equivariant measures
on $\Lcal (\GG_a,0)$ whose zeta functions lie in $\Mhat^\muboldhat\la T\ra$.
Then $\lam *\lam '$ has this property, too. If moreover 
$\lam $ and $\lam '$ have zero mass and zeta functions converging at $T=\infty$, then $\lam *\lam '(T)$ has these properties as well and 
$(\lam *\lam ')(\infty)=\lam (\infty ) *\lam'(\infty)$. 
\end{theorem}

\begin{corollary}
Let $X$ and $Y$ be smooth connected varieties and $f:X\to\GG_a$,
$g:Y\to\GG_a$ nonconstant morphisms with zero fibers $X_0$ and $Y_0$. Let $f*g
:X\times Y\to \GG_a$ be defined by $(f*g)(x,y):= f(x)+g(y)$. Then the 
restriction of  $[\phi_{f*g}]$ to $X_0\times Y_0$ and the exterior
$*$-product  $[\phi_f]*[\phi_g]\in M_{X_0\times Y_0}$ coincide. 
\end{corollary}

If we apply the Hodge number characteristic followed by formation
of the spectrum, then we recover the Thom-Sebastiani 
property for the spectrum, proved earlier by Varchenko in case $f$ and $g$ 
have isolated singularities and by M.~Saito \cite{s2} in general.

For the proof of Theorem \ref{abstractts} we need the following 

\begin{lemma}\label{residuelemma}
Let $M$ and $N$ be $\ZZ [L,L^{-1}]$-modules and let
$a\in M\la T\ra$  and $b\in N\la T\ra$ both be zero at $T=0$ and regular at 
$T=\infty$. If
$\sum_{k>0} a_kT^k$ resp.\ $\sum_{k>0} b_kT^k$ are their expansions
at $0$, then  $\sum_{k>0} (a_k\otimes b_k)T^k$ is the expansion at zero of a
$c\in (M\otimes_{\ZZ [L,L^{-1}]}N)\la T\ra$ whose value at $T=\infty$ equals 
$-a(\infty)\otimes b(\infty)$.   
\end{lemma}
\begin{proof} It is easy to see that it suffices to prove this for 
$M=N=\ZZ[ L,L^{-1}]$. 
The idea of the proof in this case is inspired by a paper of
Deligne \cite{deligne}. Fix for the moment $L\in\CC -\{0\}$.
Let $r_0>0$ be a radius of convergence for
the two expansions. Let $T\in \CC$ be such that $|T| <r_0^2$ and choose
$|T|/r_0 <r<r_0$. Consider the integral 
\[
c(T):=\frac{1}{2\pi\sqrt{-1}}\int _{|\tau |=r} a(T/\tau )b(\tau
)\frac{d\tau }{\tau }. \]
On the circle of integration the expansions converge uniformly and
absolutely and so
\[
c(T)=\frac{1}{2\pi\sqrt{-1}}\int _{|\tau |=r}\sum_{k,l \in\NN} 
a_kT^kb^l \tau ^{k-l}\frac{d\, \tau }{\tau }.
\]
Since summation and integration may be interchanged, 
only the terms with $k=l$ remain and hence $c(T)=\sum_{k\in\NN}
a_kb_kT^k$.
If $P_a$ resp.\ $P_b$ denotes the set of poles of $a$ resp.\ $b$,
then 
the integrand has polar set $TP_a^{-1} \cup P_b$ (there is no
pole in $0$ or $\infty$) and the poles enclosed by the circle of
integration are those in $TP_a^{-1}$. By the theory of residues, 
$-c(T)$ must then be equal to the sum of the residues of the
integrand at 
$P_b$. This description no longer requires $|T| <r_0^2$ and defines
an analytic extension of $c$ to the  complement of $P_aP_b$. This
extension is easily seen to be meromorphic at $P_aP_b$. 
To compute its behavior at $\infty$, we
note that $a(T/\tau )$ converges for $T\to\infty$ on a neighborhood
of $P_b$ absolutely (with all its derivatives)  to the constant function
$a(\infty)$.  So as $T\to\infty$, $-c(T)$ tends to the sum of the
residues of $a(\infty)b(\tau )\tau ^{-1}\, d\tau $ at  $P_b$. This
sum is opposite to the residue at the remaining pole $\infty$, hence equal
to $a(\infty)b(\infty )$. 
In particular, $c$ is a rational function with polar set contained in
$P_aP_b$.

Assume now that $a,b\in \tilde R$. A pole of an element of
$R$ in $\CC^\times\times\CC$ satisfies an equation $T^N=L^\nu$
for certain integers $N>0$, $\nu\ge 0$. A product of such poles
satisfies a similar equation, and this implies that a product of $c$ and a finite
set of polynomials of the form $T^N-L^\nu$ is in $\CC [L,L^{-1},T]$.
Since the expansion of $c$ at $T=0$ has integral coefficients, this product
lies in $\ZZ [L,L^{-1},T]$.
\end{proof}

\begin{proof}[Proof of Theorem \ref{abstractts}] 

We start with the convolution formula \ref{convolution}. 
It says that
\begin{multline*}
(\lam *\lam ')(T)=\\ 
-\sum_{n>0}\lam _n*\lam _nT^n +
(\LL-1)\sum_{0<i\le n} \overline{\lam _i\lam '_i}\LL^{i-n}T^n
+(\LL-1)\sum_{i>n>0} (\lam_n \lambar'_i+ \lambar_i\lam'_n)T^n. 
\end{multline*} 
We now assume that $\lam$ and $\lam'$ are massless so that
$(\LL-1)\sum_{i>n}\lambar_i=-(\LL-1) \sum_{i=1}^n\lambar_i$ 
and similarly for $\lam'$. We then have
\begin{multline*}
(\lam *\lam )(T)= -\sum_{n>0} \lam_n *\lam'_nT^n 
+(\LL-1)\sum_{0<i\le n}\overline{\lam_i\lam '_i}\LL^{i-n}T^n+\\ 
+(\LL-1)\sum_{n>0}\overline{\lam_n\lam'_n}T^n-(\LL-1)\sum_{0<i\le n}
(\lam_n\lambar'_i+\lambar_i\lam'_n)T^n.
\end{multline*} 
We consider each series on the right separately. 
By Lemma \ref{residuelemma},  $-\sum_{n>0} \lam_n *\lam'_nT^n$ is in
the in $\Mhat_k^\muboldhat\la T\ra $ with value at $\infty$ equal to 
$\lam (\infty )*\lam'(\infty )$.
We also have
\[
(\LL-1)\sum_{0<i\le n}\overline{\lam_i\lam'_i}\LL^{i-n}T^n
=(\LL-1)\sum_{i>0}\sum_{k\ge 0}
\overline{\lam_i\lam'_i}\LL^{-k}T^{k+i}=
\frac{\LL-1}{1-\LL^{-1}T}\sum_{i>0} \overline{\lam_i\lam'_i}T^i.
\]
By the same \ref{residuelemma} the righthand side is in 
$\Mhat_k^\muboldhat\la T\ra $ and takes the value zero at $\infty$.  
Since
\[
\sum_{0<i\le n} \lambar_iT^n= -(T-1)^{-1}\sum_{i>0}\lambar_iT^i
\]
is in $\Mhat_k^\muboldhat\la T\ra $ with value zero at $\infty$ it follows 
from \ref{residuelemma} that
the same is true for $(\LL-1)\sum_{0<i\le n} (\lambar_i\lam'_n)T^n$.
Likewise for $(\LL-1)\sum_{0<i\le n} (\lam_n\lambar'_i)T^n$. 
So $(\lam*\lam')(T)$ is in $\Mhat_k^\muboldhat\la T\ra $ and has value 
$\lam(\infty)*\lam'(\infty )$ at $\infty$. 
\end{proof}

\section{The McKay correspondence \cite{13}, \cite{27}, \cite{R}}

Suppose a group $G$ of finite order $m$ acts well and effectively on a
smooth connected variety $U$ of dimension $d$. This defines an orbifold
$p:U\to U_G$ with underlying variety $G\bs U$. Let us write $X$ for
the orbifold $U_G$. We also fix a primitive $m$th root of unity $\zeta_m$.

Let $g\in G$ and let $U^g$ be its fixed point set in $U$. The action of
$g$ in the normal bundle of $U^g$ decomposes that bundle into a direct sum
of eigensubbundles
\[
\nu_{U/U^g}=\oplus_{k=1}^{m-1}\nu_g^k,
\]
where $\nu_g^k$ has eigenvalue $\zeta_m^k$. We like to think of $\nu_g^k$
as the pull-back of a fractional bundle on a subvariety of $X$ whose
virtual rank is $k/m$ times that of $\nu_g^k$. A more formal discussion
involves the extension
$M_X[\LL^{1/m}]$ of $M_X$ obtained by adjoining an $m$th root of $\LL$. 
To be precise, let $w(g):=\sum_k \frac{k}{m}\rk (\nu_g^k)$,
considered as locally constant function $U^g\to\ m^{-1}\ZZ$, and let
$\LL_{U^g}^{w(g)}$ be the element of $M_{U^g}[\LL^{1/m}]\subset
M_U[\LL^{1/m}]$ that this defines. Then 
$\sum_{g\in G}\LL_{U^g}^{w(g)}$ is the image under $p^*$ of  
\[
W(X)=\sum_{[g]\in\conjclass (G)}\sum_{i\in \pi_0(U^g)}[(G_i\bs
U^g_i)/X]\LL^{w_i}
\in M_X[\LL^{1/m}].
\]
Here $U^g_i$ is the connected component of $U^g$ labeled by $i$,
$G_i$ is the $G$-stabilizer of this component, and $w_i$ the value of
$w(g)$ on $U^g_i$. The sum is over a system of representatives of the
conjugacy classes of $G$ and can be rewritten as one over  
the orbifold strata of $X$ (see Reid \cite{R}):
the decomposition of $U$ into connected strata by orbit type (a stratum is
a connected component of the locus of points with given $G$-stabilizer)
induces a partition of $X$ into orbifolds and $W(X)$ has the form $\sum_S
[S]W_S$, where the sum is over the orbifold strata, and $W_S$ is a
polynomial in $\LL^{1/m}$.
We will see that $W(X)$ can be understood as the class of an obstruction
bundle for lifting arcs in $X$ to arcs in $U$.

The McKay correspondence identifies $W(X)$ in terms of a resolution of
$X$: 

\begin{theorem}[Batyrev \cite{13}, Denef-Loeser \cite{27}] 
Let $H:Y\to X$ be a resolution of the orbifold $X$ whose
exceptional divisor $E$ has simple normal crossings. With the
usual meaning of $E^\circ_I$ and with $\nu^*_i$ as defined below we have 
the following identity in $\Mhat_X[\LL^{1/m}]$:
\[
W(X)=\sum_{I\subset \irr (D)} 
[E_I^\circ /X]\prod_{i\in I} \frac{\LL-1}{\LL^{\nu^*_i}-1}.
\]
\end{theorem}

The statement does not involve arc spaces, but the proof does.
It could well be that the identity is already valid in $M_X[\LL^{1/m}]$. 
The relative simplicity of the lefthand side has implications for the
righthand side, one of which is that all the `non-Tate' material in 
a fiber of $H$ must cancel out in the sum. For that same reason 
the lefthand side is hardly affected if we apply the weight character
relative to $X$ to it, that is, if we take the image of $W(X)$ in the
Grothendieck ring of constructible $\ZZ ((w^{-1/m}))$-modules on
$X$: just substitute $w^2$ for $\LL$.

We first seek an orbifold measure on
$\mu^\orb_{\Lcal (X)}$ on $\Lcal (X)$ with the property that for every
$G$-invariant measurable $A\subset\Lcal (U)$ we have 
\[
\mu^\orb_{\Lcal (X)} (p_*A):=\overline{\mu_{\Lcal (U)}(A)},
\]
where the righthand side should be interpreted as follows: 
think of $\mu_{\Lcal (U)}(A)$ as an element of $\Mhat_k^G$, and then let
$\overline{\mu_{\Lcal (U)}(A)}$ be the image of $\mu_{\Lcal (U)}(A)$ under
the augmentation $\Mhat_k^G\to \Mhat_k$. Since $p_*: \Lcal (U)\to \Lcal (X)$
need not be surjective, this
will not characterize the orbifold measure a priori. But it suggests
how to define it: suppose that the Jacobian ideal $\Jcal_p$ has constant
order $e$ along $A$. Then the usual measure of $\Lcal (X)$ pulled back to
$A$ is  $\LL^{-e}\mu_{\Lcal (U)}|A$. We therefore want the orbifold
measure restricted to  $p_*(A)$ to be 
the restriction of $\LL^e\mu_{\Lcal (X)}$. 
This can be done as follows.
Let $r$ be a positive integer such that $(\Omega^d_U)^{\otimes r}$ 
descends to an invertible sheaf $\omega_{X}^{(r)}$ on $X$. 
(So for every $u\in U$, $G_u$ acts
on on the tangent space $T_uU$ with determinant an $r$th root of unity.)
There is a natural homomorphism $(\Omega^d_{X})^{\otimes
r}\to\omega_X^{(r)}$ whose kernel is the torsion of 
$(\Omega^d_X)^{\otimes r}$. The image of this homomorphism has the form
$\Ical^{(r)} \omega_{X}^{(r)}$ for an ideal $\Ical^{(r)} $. We set 
\[
\mu_{\Lcal (X)}^\orb := \LL^{\ord_{\Ical^{(r)} }/r}\mu_{\Lcal (X)}.
\]
It is a measure that takes values in $\Mhat_k[\LL^{1/r}]$.

\begin{lemma}\label{orbifoldmeasure}
The pull-back of $\mu_{\Lcal (X)}^\orb$ under $p^*$ is a measure that
assigns to any $G$-invariant measurable subset $A$ of $\Lcal (U)$ the
image of $A$ under the augmentation map $\Mhat_k^G\to\Mhat_k$.
\end{lemma}
\begin{proof} If we apply $p^*$ to the identity $(\Omega^d_X)^{\otimes
r}/tors=\Ical^{(r)} \omega_{X}^{(r)}$ we
get $\Jcal_p^r(\Omega^d_U)^{\otimes r}=p^*(\Ical^{(r)}  )\omega_U^{\otimes
r}$.
Since $\Omega^d_U=\omega_U$, it follows that $p^*(\Ical^{(r)}
)=\Jcal_p^r$. So 
$\mu_{\Lcal (X)}^\orb$ pulls back under $p^*$ to
$\LL^{-\ord_{\Jcal_p}+p^*(\Ical^{(r)}  )/r}\mu_{\Lcal (U)}=\mu_{\Lcal
(U)}$.
The rest is left to the reader. 
\end{proof}

The following lemma describes the direct image of $\mu_{\Lcal (X)}^\orb$
on $X$ in terms of a resolution of $X$: let $Y\to X$ be a resolution of 
singularities with simple normal crossing divisor $E$.  We have
$H^*\omega_X^{(r)}=
\tilde\Ical^{(r)} \omega_Y^{\otimes r}$ for some fractional ideal
$\tilde\Ical^{(r)} $ on $Y$. It is known that the multiplicity $m_i$
of $E_i$ in this ideal is $>-r$. So $\nu_i^*:=1+m_i/r$ is positive.
Entirely analogous to the proof of Theorem \ref{basic} one derives:

\begin{lemma}\label{orbigusa}
The direct image of $\mu_{\Lcal (X)}^\orb$ on $X$ is represented by the
class
\[ 
\sum_{I\subset \irr (E)} 
[E^\circ_I /X]\prod_{i\in I} \frac{\LL-1}{\LL^{\nu^*_i}-1}\in
\Mhat_X[\LL^{1/r}].
\]
\end{lemma}

Let $\Lcal'(X)$ be the set of arcs in $X$ not contained in the
discriminant of $p:U\to X$. This is a subset of full measure. We decompose
$\Lcal'(X)$ according to the ramification behavior of $p:U\to X$.
Let $[m]: \DD\to \DD$ be the $m$th power map and denote the parameter of
the domain by $t^{1/m}$. We regard 
$\zeta_m$ (through its action on the domain) as generator of the Galois
group of
$[m]$.  For $\g\in\Lcal'(X)$, $\g [m]$ lifts to a morphism 
$\tilde\g:\DD\to \tilde X$ and this lift is unique up to conjugation with
$G$. Given the lift, there is a $g\in G$ such that $g
\tilde\g=\tilde\g\zeta_m$. Its conjugacy class $[g]$ in $G$
only depends on $\g$. This conjugacy class determines the isomorphism type
of the  $G$-covering over $\g$: if $m'$ is the order of $g$, then
$\g^*(p)$ is isomorphic to $G\times^{\la g\ra}\DD\to \mubold_{m'}\bs\DD$,
with $g$ acting on $\DD$ as multiplication by $\zeta^{m/m'}$.
Notice that $\tilde\g(o)$ is in the fixed point set $U^g$. The
`fractional lifts' $\tilde\g$ that so arise are like arcs 
in the total space of the normal bundle $\oplus_k\nu_g^k$ of
$U^g$ (based at the zero section) which in the $\nu_g^k$-direction
develop as $t^{k/m}$ times a power series in $t$. 

Denote the set of arcs in $\Lcal'(X)$ belonging to the conjugacy class of
$[g]$ of $g$ by $\Lcal (X,[g])$. The McKay correspondence now results
from:

\begin{lemma}
The subset $\Lcal (X,[g])$ is measurable for 
$\mu_{\Lcal (X)}^\orb$ and 
the restriction of $\mu_{\Lcal(X)}^\orb)$ to this subset
is represented by the class $[(G_g\bs U^g)/X]\LL^{w(g^{-1})}\in
M_X[\LL^{1/m}]$, where $G_g$ is the $G$-stabilizer of $U^g$. 
\end{lemma}

The proof is a calculation which we only discuss in a heuristic
fashion. The elements of $\Lcal(X,[g])$ correspond to $G_g$-orbits 
of fractional lifts as described above. In view of our definition of
orbifold measure, we need to argue that these fractional lifts are 
represented by the element $\LL_{U^g}^{w(g^{-1})}$. If $r_1,\dots
,r_{m-1}$ are positive integers, then the arcs in $\oplus_k\nu_g^k$ of
$U^g$ based at the zero section and which in the $\nu_g^k$-direction
have order $r_k$ make up a constructible subset of $\Lcal
(\oplus_k\nu_g^k)$
whose class is easily seen to be equal to $\LL_{U^g}^w$, with 
$w=\sum_k (1-r_k)\rk (\nu_g^k)$. The fact is that this also holds
for the fractional values $r_k=k/m$. So in that case we have 
$w=\sum_k (1-k/m)\rk (\nu_g^k)=w(g^{-1})$.

\section{Proof of the transformation rule \cite{25}}\label{proofs}

Let $\Xcal/\DD$ be a $\DD$-variety of pure relative dimension $d$.
The $d$th Fitting ideal of $\Omega_{\Xcal /\DD}$ defines the locus where
$\Xcal$ fails to be smooth over $\DD$; we denote that ideal by
$\Jcal (\Xcal/\DD)$. Locally this ideal is obtained as follows: if
$\Xcal$ is given as a closed subset of $(\AA^{d+l})_\DD$, then 
$\Jcal_{\Xcal/\DD}$ is the restriction to $\Xcal$ of the ideal 
generated by the determinants $\det ((\p f_j/\p x_{i_k})_{j,k=1}^l)$, 
where $f_1,\dots ,f_l$ are taken from the ideal 
$I_{\Xcal}\subset\Ocal [x_1,\dots ,x_{d+l}]$
defining $\Xcal$ and $1\le i_1<\dots <i_l\le d+l$. 

Let $\g\in\Xcal_\infty$ be such that $\g^*\Jcal (\Xcal/\DD)$ has 
finite order $e$. 
This implies that $\g$ maps $\DD^\times$ to the part $(\Xcal /\DD )_\reg$ where
$\Xcal$ is smooth over $\DD$. In particular, $\g^*\Omega_{\Xcal/\DD}$
is a $\Ocal$-module of rank $d$. Since the formation of 
a Fitting ideal commutes with base change, the $d$th Fitting ideal
of $\g^*\Omega_{\Xcal/\DD}$ will be $\mcal^e$. This means that 
the torsion of $\g^*\Omega_{\Xcal/\DD}$ has length $e$. 

It is clear that $\Der_\Ocal(\Ocal_{\Xcal,\g (o)},\Ocal)
\cong\Hom_\Ocal (\g^*\Omega_{\Xcal/\DD},\Ocal)$
is a free $\Ocal$-module of rank $d$
(where $\Ocal$ is a $\Ocal_{\Xcal,\g (o)}$-module 
via $\g*$). The fiber over $o$, 
$\Hom_\Ocal (\g^*\Omega_{\Xcal/\DD},\Ocal)\otimes_\Ocal k$,
is $d$-dimensional subspace of the Zariski tangent space 
$T_{X,\g (o)}$, which we shall denote by $\hat T_{X,\g}$.
Any $\Ocal$-homomorphism $\g^*\Omega_{\Xcal/\DD}\to\Ocal/\mcal^{n+1}$
that kills the torsion lifts to a $\Ocal$-homomorphism 
$\g^*\Omega_{\Xcal/\DD}\to\Ocal$.
This is automatic when $n\ge e$ and so $\hat T_{X,\g}$ only depends on
the $e$-jet of $\g$. The space $\hat T_{X,\g}$ has a simple geometric
interpretation: it is the `limiting position' of the tangent space along
the fibers of $\Xcal/\DD$ at the generic point of $\g (\DD)$ in the closed
point $\g (o)$. 

If $\g':\DD\to\Xcal$ has the same $n$-jet as $\g$, then 
$\g^*$ and ${\g'}^*$ differ by a homomorphism
$\Ocal_{\Xcal,\g (o)}\to\mcal^{n+1}$. The reduction modulo
$\mcal^{2(n+1)}$
of this homomorphism is a $\Ocal$-derivation, i.e., defines an element of 
$\Hom_\Ocal(\g^*\Omega_{\Xcal /\DD}, \mcal^{n+1}/\mcal^{2(n+1)})$.
Its reduction modulo $\mcal^{n+2}$ will lie in 
$\hat T_{X,\g}\otimes \mcal^{n+1}/\mcal^{n+2}$, provided 
that $n\ge e$. The next lemma shows that every element of this 
$k$-vector space so arises.

\begin{lemma}\label{noshift}
Assume that $n\ge e$. The fiber of $\pi_{n+1}\Xcal_\infty\to
\pi_n\Xcal_\infty$ over $\pi_n(\g)$ is an affine space with 
translation space $\hat T_{X,\g}\otimes_k \mcal^{n+1}/\mcal^{n+2}$. 
This defines an affine space bundle of rank $d$ over the locus of 
$\pi_n\Xcal_\infty$ defined by $\ord_{\Jcal (\Xcal/\DD)}\le n$.
\end{lemma}
\begin{proof}
Assume that $\Xcal$ is given as a closed subset of $(\AA^{d+l})_\DD$
as above. There exist $f_1,\dots ,f_l\in I_{\Xcal}$ and $1\le i_1<\dots
<i_l\le d+l$ such that the
Jacobian matrix $\det ((\p f_j/\p x_{i_k})_{j,k=1}^l)$ has order
$e$ along $\g$, whereas for any other matrix thus formed the order is $\ge
e$. By means of a coordinate change we may arrange that 
\[
\g^*df_j\equiv t^{e_j}dx_j \pmod{t^{e_j+1}(dx_{j+1},\dots ,dx_{d+l})}
,\quad j=1,\dots l, 
\]
so that $e=\sum_j e_j$. The subspace of $\AA^{d+l}_k$ spanned by the last 
$d$ basis vectors is then just $\hat T_{X,\g}$.

We investigate which $u_0\in k^{d+l}$ appear as the constant
coefficient of an $u\in \Ocal^{d+l}$ with the property that 
$\g +t^{n+1}u\in\Xcal_\infty$. We first do this for
the complete intersection  defined by  $f_1,\dots ,f_l$.
This complete intersection contains $\Xcal$ and the irreducible component
that contains the image of $\g$ lies in $\Xcal$.
So we want $f_j(\g +t^{n+1}u)=0$ for $j=1,\dots ,l$. By expanding
at $\g$ this amounts to identities of the form
\[
t^{n+1}D_\g f_j(u)+ t^{2(n+1)} F_j(u)=0, \quad  j=1,\dots ,l, 
\]
with $D_\g f_j$ the derivative of $f_j$ at $\g$ and
$F_j\in \Ocal [x_1,\dots ,x_{d+l}]$. Equivalently: 
\[
t^{-e_j}D_\g f_j(u)+ t^{n+1-e_j}F_j(u)=0,\quad  j=1,\dots ,l.   
\]
All the terms are regular and the reduction modulo $t$ yields the $j$th
unit vector in $k^{d+l}$. Hensel's lemma says that a solution
$u$ exists if and only if $u_0$ solves this set of equations modulo $t$. 
This just means that $u_0\in \hat T_{X,\g}$. 
In particular, we see that for all $k\in\NN$,
$\pi_{n+k}\pi_n^{-1}\pi_n (\g)$ is isomorphic to an affine space and hence
is irreducible. This implies that all elements of $\pi_n^{-1}\pi_n (\g)$
map to the same irreducible component of the common zero locus of
$f_1,\dots ,f_l$. It follows that $\pi_n^{-1}\pi_n
(\g)\subset\Xcal_\infty$.
The last assertion is easy. 
\end{proof}

\begin{proof}[Proof of Proposition \ref{constructible}]
Suppose that $\Xcal$ is of pure relative dimension $d$. 
Let $C_e$ denote the subset of $\Xcal_\infty$ defined
by $\ord_{\Jcal (\Xcal/\DD)}=e$. It is clear that
$C_e=\pi_e^{-1}\pi_e(C_e)$. It follows from  Greenberg's theorem
\cite{greenb2} that $\pi_e(C_e)$ is constructible. Hence $C_e$ is stable
by   Lemma \ref{noshift}. We have $\cup_e
C_e=\Xcal_\infty-(\Xcal_\sing)_\infty$. 
In view of Lemma \ref{finitecover} it now suffices
to see that $\dim \pi_e(C_e)-de\to -\infty$ as $e\to\infty$. This is not
difficult.
\end{proof} 

Let $H: \Ycal\to\Xcal$ be a $\DD$-morphism of $\DD$-varieties of pure
relative dimension $d$. Recall that the Jacobian ideal $\Jcal_H$ of $H$ is
the $0$th Fitting ideal of $\Omega_{\Ycal/\Xcal}$.
Suppose $\g\in\Ycal_\infty$ is such that $\Jcal_H$ has finite order $e$
along $\g$. Then $\g$ resp.\ $H\g$ maps the generic point $\DD^\times$ to
$(\Ycal/\DD)_\reg$ resp.\ $(\Xcal/\DD)_\reg$. We have an exact sequence of
$\Ocal$-modules
\[
(H\g)^*\Omega_{\Xcal/\DD}\to
\g^*\Omega_{\Ycal/\DD}\to \g^*\Omega_{\Ycal/\Xcal}\to 0.
\]
The base change property of Fitting ideals implies that
the length of $\g^*\Omega_{\Ycal/\Xcal}$ must be $e$. So if 
$\g^*\Omega_{\Ycal/\DD}$ is torsion free and $n\ge e$, then
the kernel of the map 
\begin{equation*}
D_\g^{(n)}:\Hom_\Ocal(\g^*\Omega_{\Ycal/\DD},\Ocal /\mcal^{n+1})\to
\Hom_\Ocal((\g H)^*\Omega_{\Xcal/\DD},\Ocal /\mcal^{n+1})
\end{equation*}
induced by the derivative of $H$ is contained in 
$\Hom_\Ocal(\g^*\Omega_{\Ycal/\DD},\mcal^{n+1-e}/\mcal^{n+1})$, 
can be identified with $\Hom_\Ocal (\g^*\Omega_{\Ycal/\Xcal},\Ocal
/\mcal^{n+1})$, and is of length $e$.  
The proof of Theorem \ref{transforule} now rests on the

\begin{keylemma}\label{shift} 
Suppose $\Ycal/\DD$ smooth and let $A\subset \Ycal_\infty$ be a stable
subset of level $l$: $A=\pi_l^{-1}\pi_l(A)$. Assume that $H\big|_A$ is
injective and that $\ord_{\Jcal_H}\big|_A$ is constant equal to $e<\infty
$.
If $n\ge \sup\{2e,l+e,\ord_{\Jcal (\Xcal/\DD)}\big|_{HA} \}$, then 
$H_n: \pi_nA\to H_n\pi_nA$
has the  structure of affine-linear bundle of dimension $e$.
\end{keylemma}
\begin{proof}
Let $\g\in A$ and put $x:=\g (o)$, $y:=H(x)$. Suppose $\g'\in A$ is such
that $H\g'$ and $H\g $ have the same $n$-jet.
We first show that $\g$ and $\g'$ have the same $(n-e)$-jet. We do this
by constructing a $\g_1\in\Ycal_\infty$ (by successive approximation) with
the same $(n-e)$-jet as $\g$ and with $H\g_1=H\g'$.
Since $n-e\ge l$, we will have $\g_1\in A$ and our injectivity assumption
then implies $\g_1=\g'$.

The difference $(H\g)^*-(H\g')^*$ defines a
$\Ocal$-derivation $\Ocal_{\Xcal,x}\to\mcal^{n+1}/\mcal^{2(n+1)}$ over
$\g^*$ and hence a 
$\tilde v\in\Hom_\Ocal ((\g H)^*\Omega_{\Xcal
/\DD},\mcal^{n+1}/\mcal^{2n+2})$. Since $n\ge \ord_{H\g}\Jcal
(\Xcal/\DD)$, 
this element annihilates the torsion of  $(\g H)^*\Omega_{\Xcal /\DD}$. 
This is then also true for its reduction modulo $\mcal^{n+2}$ and it 
follows from the fact that $n\ge e$ that this reduction is of the form
$D_\g^{(n+1)}(u)$ for some 
$u\in\Hom_\Ocal (\g^*\Omega_{\Ycal/\DD},\mcal^{n-e+1}/\mcal^{n+2})$. 
Regard $u$ as a $\Ocal$-derivation 
$\Ocal_{\Ycal,y}\to\mcal^{n-e+1}/\mcal^{n+2}$ and let
$\g_1\in\Ycal_\infty$ 
be such that $\g_1^* -\g^*$ represents $u$. Then $\pi_{n-e}(\g_1)=
\pi_{n-e}(\g)$ and $\pi_{n+1}(H\g_1)=\pi_{n+1}(H\g)$.
Replace $\g$ by $\g_1$ and continue with induction on $n$.

So $(\g')^*-\g^*$ defines a
$\Ocal$-derivation  $\Ocal_{\Xcal,x}\to \mcal^{n-e+1}/\mcal^{2(n-e+1)}$
and hence a
$\Ocal$-derivation $\Ocal_{\Xcal,x}\to \mcal^{n-e+1}/\mcal^{n+1}$ 
(because $n\ge 2e$). The latter is zero if and only if
$\pi_n(\g')=\pi_n(\g)$. 
This proves that the fiber of $H_n\big|_{\pi_nA}$ through $\pi_n(\g)$
is an affine space over the kernel of $D_\g^{(n)}$, 
$\Hom_\Ocal (\g^*\Omega_{\Ycal/\Xcal},\Ocal/\mcal^{n+1})$, which has
length $e$.

The last assertion is easy.
\end{proof}

\begin{proof}[Proof of \ref{transforule}] 
It is enough to prove this for $A$ stable.
In that case the theorem follows in a straightforward manner from
Lemma \ref{shift}.
\end{proof}

\end{document}